\documentclass[11pt]{amsart}
\usepackage{amsmath,amssymb,amsthm,mathtools}
\usepackage[margin=1.15in]{geometry}
\usepackage{enumitem}
\usepackage{stmaryrd}
\usepackage[colorlinks=true,linkcolor=blue,citecolor=blue,urlcolor=blue]{hyperref}

\newtheorem{theorem}{Theorem}
\newtheorem{thmx}{Theorem}
\newtheorem{proposition}[theorem]{Proposition}
\newtheorem{lemma}[theorem]{Lemma}
\newtheorem{corollary}[theorem]{Corollary}
\theoremstyle{definition}
\newtheorem{definition}[theorem]{Definition}
\theoremstyle{remark}
\newtheorem{remark}[theorem]{Remark}
\numberwithin{equation}{section}
\numberwithin{theorem}{section}

\newcommand{\R}{\mathbb{R}}
\newcommand{\Sph}{\mathbb{S}}
\newcommand{\RP}{\mathbb{RP}}
\newcommand{\Z}{\mathbb{Z}}
\newcommand{\Ric}{\operatorname{Ric}}
\newcommand{\Hess}{\operatorname{Hess}}
\newcommand{\Riem}{\operatorname{Riem}}
\newcommand{\dv}{\operatorname{div}}
\newcommand{\inte}{\operatorname{int}}

\title[Area-minimizing surfaces and electrostatic systems]{Area-minimizing surfaces and the rigidity\\ of three-dimensional electrostatic systems}
\author{Rafael Rocha de Farias}
\address{Universidade Federal do Par\'a, Instituto de Ci\^encias Exatas e da Natureza, Bel\'em, Par\'a, 66075-110, Brazil}
\email{rafaelfariasmat@ufpa.br}

\begin{document}
\begin{abstract}
We study three-dimensional electrostatic systems, the Riemannian data of static solutions of the Einstein--Maxwell equations with cosmological constant $\Lambda$, by means of locally area-minimizing surfaces. Such a surface with non-empty boundary on the horizons is a totally geodesic hemisphere and forces the de Sitter system. A closed one in the interior forces a system of product type, one of eight models if the scalar curvature is somewhere positive on it, among them the Bertotti--Robinson, Nariai and cold and ultracold black hole systems. No compactness, asymptotic or pinching assumption is needed. In particular, closed systems without horizons are flat or compact quotients of the Nariai--Bertotti--Robinson system. For compact systems with connected horizon, the electric field vanishes even when the system is not simply connected, and the area of the horizon is either $12\pi/\Lambda$, for de Sitter, or at most $4\pi/\Lambda$. With several horizons, every horizon whose area does not exceed the sum of the others is stable, minimizes area in its homology class and satisfies an area--charge inequality, with equality only for the charged Nariai systems; its charged Hawking mass and charge therefore lie in the parameter region of the Reissner--Nordstr\"om--de Sitter black holes. Compact systems with horizons are spherical space forms with balls removed. Noncompact systems with strictly mean-convex ends have one end, and their horizons are stable and $\pi_1$-injective and satisfy a strict area--charge inequality.
\end{abstract}
\maketitle

\section{Introduction}\label{sec:intro}

A \emph{$3$-dimensional electrostatic system} with cosmological constant $\Lambda\in\R$ is a quadruple $(M^3,g,V,E)$, where $(M,g)$ is a connected, complete, orientable Riemannian $3$-manifold with or without boundary, $V\in C^\infty(M)$ is positive in $\inte M$, and $E\in\mathfrak{X}(M)$ is a smooth vector field, such that
\begin{equation}\label{eq:electrostatic}
\left\{\begin{aligned}
\Hess_gV&=V\bigl({\Ric_g}-\Lambda g+2E^\flat\otimes E^\flat-|E|_g^2\,g\bigr),\\
\Delta_gV&=\bigl(|E|_g^2-\Lambda\bigr)V,\\
\dv_gE&=0,\qquad d(VE^\flat)=0,
\end{aligned}\right.
\end{equation}
and $V^{-1}(0)=\partial M$; in particular, $V>0$ on $M$ if $\partial M=\emptyset$. These are the equations for which $(\R\times\inte M,-V^2dt^2+g)$, with Faraday tensor $F=V\,E^\flat\wedge dt$, is a static solution of the Einstein--Maxwell equations with cosmological constant $\Lambda$ \cite{CLS24}; the components of $\partial M$, which we call \emph{horizons}, are the bifurcation surfaces of its Killing horizons. Tracing the first equation and using the second gives $R_g=2(\Lambda+|E|^2)$. The equations are linear in $V$, so $(M,g,cV,E)$ is again an electrostatic system for every constant $c>0$ (Lemma~\ref{lem:basic}(iii)), with the same spacetime after the rescaling $t\mapsto t/c$ of time. We do not distinguish between $(M,g,V,E)$ and $(M,g,cV,E)$; statements such as ``the system is the de Sitter system'' are understood in this sense. In \cite[Definition~1]{CLS24} the potential may change sign; our systems are the closures of the regions where it is positive, so we recall the condition $V^{-1}(0)=\partial M$ in the statements, and systems on which $V$ vanishes in the interior are treated in the companion paper \cite{Far26b}.

Besides the vacuum ($E=0$) de Sitter, Nariai \cite{Nar51} and Schwarzschild--de Sitter systems, the class contains the Reissner--Nordstr\"om--de Sitter, charged Nariai, cold and ultracold black hole systems, and the Bertotti--Robinson-type cylinders \cite{Ber59,Rob59}. For $\Lambda=0$ and an asymptotically flat end, such systems are classified, under suitable assumptions on the charges, by the uniqueness theorems for static electrovacuum black holes; see \cite{CT07} and the references there. For $\Lambda>0$ one considers compact systems with horizons, and a central problem is the static form of the cosmic no-hair conjecture, due to Boucher, Gibbons and Horowitz \cite{BG83,BGH84}: \emph{the de Sitter system is the only compact simply connected static vacuum system with positive scalar curvature and connected boundary}. The sharp estimate $\Lambda|\partial M|\le12\pi$ is known \cite{BGH84,Shen97}, as are several partial results \cite{Kob82,Laf83,Amb17}; to our knowledge, the conjecture is open \cite{FLRS25}. Freitas, Leandro, Ribeiro and Sabo investigated whether the de Sitter system is the unique compact simply connected electrostatic system with positive scalar curvature \cite[Section~1]{FLRS25}, an electrostatic analogue of this conjecture. For a connected horizon, this question is equivalent to the conjecture itself, by Gauss's law: if $M$ is simply connected, then $VE^\flat=d\psi$ with $\psi$ constant on the horizon, the flux of $E$ through the horizon vanishes, and hence $\int_MV|E|^2=0$ (compare Lemma~\ref{lem:potential}). Without that assumption the argument is not available, and the results obtained so far use pinching conditions on $|E|$ \cite[Theorem~E]{CLS24}, \cite[Theorems~A, C, D and~F]{FLRS25}. Theorem~\ref{thm:C} below shows that $E\equiv0$ on every compact system with connected horizon, simply connected or not. Such a system is therefore a static vacuum system, and for connected horizons the pinching conditions above always hold.

Since the work of Schoen and Yau \cite{SY79}, locally area-minimizing surfaces have been used to prove rigidity results under lower bounds on the scalar curvature \cite{CG00,BBN10,Nun13}, also for charged initial data sets \cite{BBB23,GM25,CM25,Men25}. For static systems, Ambrozio \cite{Amb17} and Cruz, Lima and Sousa \cite{CLS24} observed that the potential provides a canonical deformation: flowing a minimal surface in the normal direction with speed $V$ does not increase area. We state the result of \cite{CLS24} in the form used in this paper: compactness of $M$ is not needed (Lemma~\ref{lem:flow}), and in~(a), \cite{CLS24} expresses the curvature in a different form.

\begin{theorem}[{\cite[Proposition~7]{CLS24}}; see {\cite[Proposition~14]{Amb17}} for $E=0$]\label{thm:CLS}
Let $(M,g,V,E)$ be an electrostatic system with $V^{-1}(0)=\partial M$, and let $\Sigma_0\subset M$ be a compact, connected, orientable minimal surface which is either closed and contained in $\inte M$, or almost properly embedded with $\Sigma_0\cap\partial M=\partial\Sigma_0\neq\emptyset$ (Definition~\ref{def:almost}). Let $\Phi:[0,\delta)\times\Sigma_0\to M$ solve $\partial_t\Phi=(V\circ\Phi)\,N_\Phi$, where $N_\Phi$ is the unit normal of $\Sigma_t=\Phi(t,\Sigma_0)$. Then $t\mapsto|\Sigma_t|$ is non-increasing. If it is constant, then $|E|$ is constant on $U=\Phi([0,\delta)\times\Sigma_0)$, each $\Sigma_t$ is totally geodesic, $E$ is normal to the leaves, and $\Phi$ is an isometry from $([0,\delta)\times\Sigma_0,(V\circ\Phi)^2dt^2+g_{\Sigma_0})$ (with $\partial\Sigma_0$ removed in case (b)) onto $(U,g)$. Moreover:
\begin{enumerate}
\item[(a)] if $\Sigma_0$ is closed, $V$ is constant on each $\Sigma_t$, so that $g|_U$ is a Riemannian product, and $(\Sigma_t,g_{\Sigma_t})$ has constant Gaussian curvature (equal to $\Lambda+|E|^2$, see Lemma~\ref{lem:product});
\item[(b)] if $\partial\Sigma_0\neq\emptyset$, then $\Lambda>0$, $E|_U=0$, $(\Sigma_t,g_{\Sigma_t})$ is isometric to a domain bounded by geodesic arcs in the round sphere of curvature $\Lambda/3$, and $V_t=V|_{\Sigma_t}$ is a static potential on $\Sigma_t$, i.e.\ $\Hess_{\Sigma_t}V_t=-\tfrac{\Lambda}{3}V_t\,g_{\Sigma_t}$.
\end{enumerate}
\end{theorem}

In~(b), \cite{CLS24} also asserts that $g|_U$ has constant sectional curvature $\Lambda/3$; this uses that $V_t$ is a multiple of $V_0$, which follows once $\Sigma_0$ is known to be a hemisphere (Lemma~\ref{lem:surfhemisphere}), and we derive it in the proof of Theorem~\ref{thm:A}. In this paper we derive global consequences of this local statement, for every $\Lambda\in\R$ and without compactness, asymptotic or pinching assumptions. Locally area-minimizing surfaces are understood as in Definition~\ref{def:lam}, and the systems mentioned below are described in Section~\ref{sec:models}.

\subsection{Rigidity and classification}

\begin{thmx}[Rigidity of locally area-minimizing surfaces with boundary]\label{thm:A}
Let $(M,g,V,E)$ be an electrostatic system with $V^{-1}(0)=\partial M$, and let $\Sigma\subset M$ be a compact, connected, orientable surface with $\Sigma\cap\partial M=\partial\Sigma\neq\emptyset$ which is almost properly embedded and locally area-minimizing (Definitions~\ref{def:almost} and~\ref{def:lam}). Then $(M,g,V,E)$ is the de Sitter system.
\end{thmx}

Theorem~\ref{thm:A} is a global rigidity statement whose only hypothesis is variational: compactness of $M$, the sign of $\Lambda$ and the vanishing of $E$ are consequences. Its proof shows that $\Sigma$ is a totally geodesic hemisphere bounded by a great circle of $\partial M$, a page of the open book foliation of Remark~\ref{rem:openbook}; in particular, $\Sigma$ is a disk with $\Lambda|\Sigma|=6\pi$, and no annulus satisfies the hypothesis. The hypothesis holds for the least area disks and annuli given by the Plateau and Douglas problems (Theorems~\ref{thm:plateau} and~\ref{thm:douglas} and Remark~\ref{rem:almostflow}), so Theorem~\ref{thm:A} bounds the area of such disks and excludes such annuli (Theorem~\ref{thm:F}(iv) and Corollary~\ref{cor:disks}).

\begin{thmx}[Classification]\label{thm:B}
Let $(M,g,V,E)$ be an electrostatic system with $V^{-1}(0)=\partial M$, possibly $\partial M=V^{-1}(0)=\emptyset$, such that $\inte M$ contains a closed embedded surface $\Sigma$ which is either orientable and locally area-minimizing, or non-orientable with locally area-minimizing orientable preimage in the double cover of $M$ associated with the mod~$2$ intersection number with $\Sigma$ (Definition~\ref{def:lam}). Then $(M,g,V,E)$ is isometric to one of the systems of product type of Section~\ref{sec:models}, in all of which $g=ds^2+h$ and $E=E_0\partial_s$, with $E_0$ constant and $h$ of constant curvature $K=\Lambda+E_0^2$: one of the eight systems (i)--(viii) of Table~\ref{tab:spherical} if $R_g>0$ somewhere on $\Sigma$, which is automatic if $\Lambda>0$, and one of the systems of Table~\ref{tab:flathyp} otherwise.
\end{thmx}

Unlike the uniqueness theorems for $\Lambda=0$, Theorem~\ref{thm:B} needs no asymptotic assumption, and, unlike \cite{FLRS25}, no pinching of $E$. For $\Lambda=0$, it shows that the system is one of the Bertotti--Robinson systems (i)--(iii) if $E\not\equiv0$, and one of the flat systems (iv$'$) and (vi$'$) if $E\equiv0$. For static vacuum metrics with $\Lambda=0$, Huang, Martin and Miao \cite[Proposition~5]{HMM18} showed with the same flow that a closed locally area-minimizing hypersurface on which $V>0$ has a Ricci-flat neighborhood, and Coutinho and Leandro \cite{CoLe22} extended their argument to electrostatic systems with $\Lambda=0$, under assumptions on the scalar curvature along the surface which exclude the Bertotti--Robinson systems. For general sub-static systems, Colombo, Freitas, Mari and Rigoli \cite[Theorem~B]{CFMR26} proved that such a hypersurface has a neighborhood isometric to a warped product, and that an end relative to it splits globally if it is complete for $V^{-2}g$ and $V^2g$; Theorem~\ref{thm:B} needs no condition on the ends (Remark~\ref{rem:substatic}). The systems of Table~\ref{tab:flathyp} occur only for $\Lambda\le0$, and none of them is compact with non-empty boundary. The associated spacetimes are, up to the quotients, the products $AdS_2\times\Sph^2$, $\R^{1,1}\times\Sph^2$ and $dS_2\times\Sph^2$ for Table~\ref{tab:spherical} (Remark~\ref{rem:nearhorizon}), and $AdS_2\times\Sigma$ and $\R^{1,1}\times T^2$ for Table~\ref{tab:flathyp}.

\begin{remark}[Constant area along the flow]\label{rem:constantarea}
Theorems~\ref{thm:A} and~\ref{thm:B} use the minimizing property of $\Sigma$ only to show that $\Sigma$ has constant area along the flow: $t\mapsto|\Sigma_t|$ is constant on some interval $[0,\delta)$ for the flow of Theorem~\ref{thm:CLS} in the direction of one unit normal of $\Sigma$ (Definition~\ref{def:constantarea}). By Theorem~\ref{thm:CLS}, it suffices that $|\Sigma_t|\ge|\Sigma|$ for small $t\ge0$, a condition on a single family of competitors on one side of $\Sigma$. Both theorems hold under this weaker hypothesis, and Sections~\ref{sec:flow} and~\ref{sec:AB} prove them in this form: Theorem~\ref{thm:A} for compact, connected, orientable minimal surfaces with non-empty boundary which are almost properly embedded, with bounded second fundamental form away from the exceptional set, and Theorem~\ref{thm:B} for closed embedded minimal surfaces in $\inte M$ which are orientable or have orientable preimage with this property. For instance, a component of a homological minimizer lying on the boundary of the region of minimization has constant area along the flow (Lemma~\ref{lem:noproduct}).
\end{remark}

\subsection{Several horizons}

\begin{thmx}[Compact systems with several horizons]\label{thm:D}
Let $(M,g,V,E)$ be a compact electrostatic system with $V^{-1}(0)=\partial M=\partial_1M\cup\dots\cup\partial_\ell M$, $\ell\ge1$. Then $\Lambda>\min_M|E|^2\ge0$ and $R_g>0$ (Lemma~\ref{lem:basic}(v)), and:
\begin{enumerate}
\item[(i)] each $\partial_iM$ is homeomorphic to $\Sph^2$ \cite[Theorem~A]{CLS24}, $H_1(M;\R)=0$, and $H_2(M;\R)\cong\R^{\ell-1}$;
\item[(ii)] $VE^\flat=d\psi$ for some $\psi\in C^\infty(M)$ which is constant on each horizon;
\item[(iii)] if $\partial_JM$ is a non-empty union of horizons with $|\partial_JM|\le|\partial M\setminus\partial_JM|$, every closed embedded orientable surface homologous to $\partial_JM$ has area at least $|\partial_JM|$, and, unless $(M,g,V,E)$ is a charged Nariai system, equality holds only for $\partial_JM$ and, if $|\partial_JM|=|\partial M\setminus\partial_JM|$, for $\partial M\setminus\partial_JM$. In particular, every horizon of area at most the sum of the areas of the others is stable, so an unstable horizon has larger area than all the others together;
\item[(iv)] at most one horizon is unstable; if one is, $M$ is homeomorphic to a $3$-ball minus $\ell-1$ open $3$-balls \cite[Theorems~A and~10]{CLS24}, and it has Morse index one;
\item[(v)] $M$ is diffeomorphic to a spherical space form $\Sph^3/G$ with $\ell$ open balls removed; in particular $\pi_1(M)=G$ is finite, and $G$ is trivial if some horizon is unstable.
\end{enumerate}
\end{thmx}

Item (iii) says that unions of horizons carrying at most half of the total area are outer-minimizing, the property assumed in Penrose-type inequalities \cite{HI01,Bra01}; for black holes inside an unstable cosmological horizon, every surface enclosing some of them has area larger than the sum of their areas. It is also a quantitative form of the uniqueness of the unstable horizon in \cite[Theorem~B]{Amb17} and \cite[Theorem~A]{CLS24}. For simply connected static vacuum systems, Ambrozio's inequality $\sum_i\kappa_i(\Lambda|\partial_iM|-4\pi)<0$ \cite[Theorem~C]{Amb17}, valid unless the system is the de Sitter or the Nariai system, bounds an average of the areas weighted by the surface gravities, whereas (iii) and Theorem~\ref{thm:masscharge} bound each horizon except one of largest area, and every horizon if $M$ is not simply connected. The index-one statement in (iv) is derived from Theorem~\ref{thm:A}, unlike the min-max results of \cite[Theorems~B and~C]{CLS24} (Remark~\ref{rem:CLSindex}), and (v) uses the topological classification of Carlotto and Li \cite{CL24}. When $\pi_1(M)\neq1$, every sphere separating the summand $\Sph^3/G$ from all horizons has area at least $|\partial M|$ (Corollary~\ref{cor:spheres}).

For a horizon $S$ with charge $Q_S=\frac1{4\pi}\int_S\langle E,\nu\rangle$ and area radius $r_S=(|S|/4\pi)^{1/2}$, the \emph{charged Hawking mass} of the minimal surface $S$ is
\[
m_H(S)=\frac{r_S}{2}\Bigl(1+\frac{Q_S^2}{r_S^2}-\frac{\Lambda r_S^2}{3}\Bigr);
\]
it depends only on $|S|$ and $Q_S$, and it equals the mass parameter $m$ at every horizon of a Reissner--Nordstr\"om--de Sitter system. For $\Lambda>0$ and $4\Lambda Q^2\le1$ let
\[
r_\pm(\Lambda,Q)=\Bigl(\frac{1\pm\sqrt{1-4\Lambda Q^2}}{2\Lambda}\Bigr)^{1/2},\qquad m_\pm(\Lambda,Q)=r_\pm\Bigl(1-\frac23\Lambda r_\pm^2\Bigr),
\]
so that $m_-(\Lambda,0)=0$ and $m_+(\Lambda,0)=1/(3\sqrt\Lambda)$. For fixed $\Lambda$ and $Q$, the Reissner--Nordstr\"om--de Sitter metrics of mass $m$ have a black-hole horizon and a cosmological horizon exactly when $m_-<m<m_+$; the limits $m_-$ and $m_+$ correspond to the cold black hole and to the charged Nariai system \cite[Section~3.4]{CLS24}, and $r_\pm$ are the critical radii at which the stability of the horizons of these metrics changes \cite[Section~3.8]{CLS24}.

\begin{thmx}[Area, mass and charge of horizons]\label{thm:masscharge}
Let $(M,g,V,E)$ be a compact electrostatic system with $V^{-1}(0)=\partial M$ consisting of $\ell\ge2$ horizons, and let $S$ be a horizon with $|S|\le|\partial M|-|S|$; this holds for every horizon except possibly one of largest area. Then $S$ is stable, $4\Lambda Q_S^2<1$ and
\[
\Lambda|S|+\frac{16\pi^2Q_S^2}{|S|}\le4\pi ,
\]
with equality if and only if $(M,g,V,E)$ is a charged Nariai system. In terms of the area radius and the charged Hawking mass, this means that
\[
r_-(\Lambda,Q_S)<r_S\le r_+(\Lambda,Q_S)\qquad\text{and}\qquad m_-(\Lambda,Q_S)<m_H(S)\le m_+(\Lambda,Q_S),
\]
with equality in either upper bound if and only if $(M,g,V,E)$ is a charged Nariai system. Moreover, the charges $Q_j=Q_{\partial_jM}$ satisfy $\sum_jQ_j=0$, at least two of them are non-zero if $E\not\equiv0$, and $|Q_j|<(\ell-1)/(2\sqrt\Lambda)$ for every $j$. If $M$ is not simply connected, these bounds hold, with strict inequalities, for every horizon.
\end{thmx}

Thus every horizon of a compact system, except possibly one of largest area when $M$ is simply connected, is stable, and its charged Hawking mass and charge lie in the parameter region of the Reissner--Nordstr\"om--de Sitter black holes, although no symmetry is assumed.

For two horizons, as in the Schwarzschild--de Sitter and Reissner--Nordstr\"om--de Sitter systems, where $\partial_1M$ is the black-hole horizon and $\partial_2M$ the cosmological one, labelled so that $|\partial_1M|\le|\partial_2M|$, the hypothesis of Theorem~\ref{thm:masscharge} holds for $\partial_1M$. Thus $\partial_1M$ is always stable and satisfies the area--charge inequality, with equality only for the charged Nariai systems. An unstable horizon can only be $\partial_2M$, and then $|\partial_1M|<|\partial_2M|$ and $M$ is diffeomorphic to $\Sph^2\times[0,1]$ (Theorem~\ref{thm:D}(iii)--(v)), and its area and charge are bounded in Remark~\ref{rem:largest}; in particular, horizons of equal area are both stable. Moreover $Q_1=-Q_2$ and $\int_MV|E|^2=-4\pi(\psi_1-\psi_2)Q_1$ (Lemma~\ref{lem:potential}), so that $E\not\equiv0$ if and only if the potential $\psi$ of Theorem~\ref{thm:D}(ii) takes different values $\psi_1,\psi_2$ on the two horizons, and then $Q_1\neq0$ and $(\psi_1-\psi_2)Q_1<0$.

The area--charge inequality holds for every stable horizon, as for stable marginally trapped surfaces \cite{DJR12,Sim12}, \cite[Proposition~3.1]{GM25}, and for stable minimal surfaces it is equivalent to the lower bound for the charged Hawking mass of Leandro and Sabo \cite{LeSa25}. For two horizons, the bounds $4\Lambda Q_S^2\le1$ and $r_S\ge r_-$ for the smaller horizon $S$, with equality only for a product, also follow from \cite[Theorem~1]{CM25}, which holds for initial data with $R_g\ge2\Lambda+2|E|^2$ and assumes $H_2(M,S)=0$. By Lefschetz duality, $H_2(M,S)\cong H^1(M,\partial M\setminus S)$, and since $H^1(M)=0$ by Theorem~\ref{thm:D}(i), this group vanishes if and only if $\partial M\setminus S$ is connected, that is, $\ell=2$; for $\ell\ge3$ that hypothesis fails. What is new is that the condition $|S|\le|\partial M|-|S|$ alone makes $S$ area-minimizing (Theorem~\ref{thm:D}(iii)), which gives its stability, the upper bounds, the strict inequalities and the characterization of the charged Nariai systems, together with the constraints on the charges. The equality case uses the local rigidity of area-minimizing horizons of Cruz and Mendes \cite[Proposition~6]{CM25} (Proposition~\ref{prop:rigidity}), a charged version of \cite{BBN10} (see also \cite{BBB23,Men25}). Once $S$ is known to minimize area, the upper bound and the characterization also follow from the global splitting theorem of Galloway and Mendes \cite[Theorem~4.1]{GM25}, which holds for initial data with $R_g\ge2\Lambda+2|E|^2$ and $\Lambda>0$, whenever $H_2(M,S)=0$ and $M$ satisfies their homotopy condition with respect to $S$, for instance when $\ell=2$ and $M\cong\Sph^2\times[0,1]$.

\subsection{Noncompact systems}

We say that $(M,g)$ has \emph{mean-convex ends} if $M$ is compact or admits an exhaustion by compact domains $\Omega_k\supset\partial M$ (in particular, $\partial M$ is compact) whose outer boundaries $\partial_+\Omega_k=\partial\Omega_k\setminus\partial M$ are smooth, with mean curvature vector pointing into $\Omega_k$ or vanishing; the ends are \emph{strictly mean-convex} if this vector never vanishes, as for asymptotically Euclidean or hyperbolic ends. A system is \emph{of product type} if it is a quotient of a product system (Definition~\ref{def:product}); the models of Theorem~\ref{thm:B} are of this type. We say that $\partial M$ is \emph{strictly outer-minimizing} if $|\Theta|<|\tilde\Sigma|$ for every non-empty union $\Theta$ of components of $\partial M$ and every closed embedded orientable surface $\tilde\Sigma\neq\Theta$ homologous to $\Theta$, as in the Riemannian Penrose inequality \cite{HI01,Bra01}. Two disjoint smooth Jordan curves $\Gamma_1,\Gamma_2$ satisfy the \emph{strict Douglas condition} if $\inf_A|A|<\inf_{D_1,D_2}(|D_1|+|D_2|)$, over immersed annuli $A$ with $\partial A=\Gamma_1\cup\Gamma_2$ and immersed disks $D_i$ with $\partial D_i=\Gamma_i$; this is the classical hypothesis for the existence of least area annuli \cite{Dou31,Jos85}, \cite[Theorem~1]{MY82a}.

\begin{thmx}[Noncompact systems]\label{thm:F}
Let $(M,g,V,E)$ be a noncompact electrostatic system with $V^{-1}(0)=\partial M$ and strictly mean-convex ends.
\begin{enumerate}
\item[(i)] If $\partial M\neq\emptyset$, then every component of $\partial M$ is a stable minimal surface and $\partial M$ is strictly outer-minimizing.
\item[(ii)] $M$ has exactly one end, and $H_2(M;\Z)$ is generated by the classes of the components of $\partial M$; in particular, $H_2(M;\Z)=0$ if $\partial M=\emptyset$. Every closed embedded orientable surface in $\inte M$ separates $M$, and the flux of $E$ through it is $4\pi\sum_jn_jQ_j$ for some integers $n_j$, where $Q_j$ is the charge of $\partial_jM$.
\item[(iii)] The inclusion of every component of $\partial M$ is injective on $\pi_1$; in particular, if $\pi_1(M)$ contains no fundamental group of a closed orientable surface of positive genus, the components of $\partial M$ are spheres.
\item[(iv)] If $\Gamma_1,\Gamma_2$ are smooth Jordan curves in different components of $\partial M$, null-homotopic in $M$, then no annulus bounded by $\Gamma_1\cup\Gamma_2$ has least area among such annuli; in particular, the strict Douglas condition fails for $\Gamma_1,\Gamma_2$.
\end{enumerate}
Item (iv) also holds when $M$ is compact.
\end{thmx}

Strict mean-convexity is needed in (i) and (ii): the cold and ultracold black holes and the cylinders of Theorem~\ref{thm:B} have mean-convex ends, but not strictly mean-convex ones; the horizons of the former are not strictly outer-minimizing, and the latter have two ends (Remark~\ref{rem:sharp}). For $\Lambda>0$ we know no noncompact system with strictly mean-convex ends, so (i) and (ii) concern mainly $\Lambda\le0$, for instance the Schwarzschild and subextremal Reissner--Nordstr\"om exteriors and their anti-de Sitter analogues. By (ii), and in contrast with the compact quotients of the Nariai--Bertotti--Robinson system in Theorem~\ref{thm:E}(b), whose slices $\{s\}\times\Sph^2$ carry a non-zero flux, the electric flux through every closed surface is determined by the charges of the horizons. Items (ii) and (iii) are Riemannian, static analogues of topological censorship and of Hawking's theorem on the topology of horizons \cite{Haw72,GSWW99}; (iii) allows horizons of positive genus, as in the topological black holes with $\Lambda<0$ \cite{GSWW99}, but forces them to be $\pi_1$-injective. Moreover, every horizon $S$ in (i) satisfies the strict area--charge inequality $\Lambda|S|+16\pi^2Q_S^2/|S|<2\pi\chi(S)$; hence the horizons are spheres if $\Lambda\ge0$, and $|S|>4\pi Q_S^2$ if $\Lambda=0$ (Corollary~\ref{cor:noncompactcharge}; compare \cite[Theorem~3]{CM25}). Finally, in systems with mean-convex ends, compact or not, a disk spanning a curve on a horizon has at least the least area of a disk bounded by that curve in the horizon, and an annulus joining curves on two horizons has more than the sum of the corresponding areas, which is the infimum (Corollary~\ref{cor:disks}).

\subsection{Closed systems}

\begin{thmx}[Closed systems]\label{thm:E}
Let $(M,g,V,E)$ be a compact electrostatic system with $\partial M=V^{-1}(0)=\emptyset$. Then either
\begin{enumerate}
\item[(a)] $\Lambda=0$, $E\equiv0$, $V$ is constant, and $(M,g)$ is a flat manifold, i.e.\ a finite quotient of a flat $3$-torus; or
\item[(b)] $\Lambda>0$ and $(M,g,V,E)$ is one of the compact quotients (v) of the Nariai--Bertotti--Robinson system (Section~\ref{sec:models}).
\end{enumerate}
\end{thmx}

In (b), $E\not\equiv0$, since otherwise integrating $\Delta V=-\Lambda V$ would give $\Lambda\int_MV=0$. Theorem~\ref{thm:E} determines the static solutions of the Einstein--Maxwell equations with cosmological constant and purely electric field whose spatial slices are closed and contain no horizon. In vacuum, the integration above forces $\Lambda=0$ and a flat spacetime. An electric field adds only the compact quotients of the Pleba\'nski--Hacyan spacetime $\R^{1,1}\times\Sph^2$ \cite{PH79,GP09}, with slices $\Sph^1\times\Sph^2$ and $\Lambda=E_0^2>0$ (Section~\ref{sec:related}). In both cases $V$ is constant, so that the spacetime is ultrastatic.

If $V$ is allowed to vanish in the interior, as in \cite[Definition~1]{CLS24}, the closures of the components of $\{V\neq0\}$ are compact systems with horizons, to which Theorems~\ref{thm:D} and~\ref{thm:C} apply. In the companion paper \cite{Far26b} we use this to describe the topology of closed manifolds on which $V$ changes sign and the charges of the components of $V^{-1}(0)$, and we show that adjacent regions are mirror images of each other.

\subsection{Connected horizon}

\begin{thmx}[Connected horizon]\label{thm:C}
Let $(M,g,V,E)$ be a compact electrostatic system whose horizon $\partial M=V^{-1}(0)$ is non-empty and connected. Then $\partial M\cong\Sph^2$, $H_1(M;\R)=0$, $\Lambda>0$, and $E\equiv0$. In particular $(M,g,V)$ is a compact oriented static vacuum system with $R_g=2\Lambda>0$ and connected boundary. Moreover, exactly one of the following holds:
\begin{enumerate}
\item[(a)] $\partial M$ is unstable. Then $M$ is diffeomorphic to a $3$-ball, $\partial M$ has Morse index one, and either the system is the de Sitter system or $\Lambda|\partial M|<4\pi$.
\item[(b)] $\partial M$ is stable and $M$ is simply connected. Then $M$ is diffeomorphic to a $3$-ball and $\Lambda|\partial M|<4\pi$.
\item[(c)] $M$ is not simply connected. Then $\partial M$ is stable, $M$ is diffeomorphic to a spherical space form $\Sph^3/G$, $G\neq1$, with a ball removed, $\Lambda|\partial M|\le4\pi$, and every embedded projective plane $P\subset\inte M$ satisfies $|P|\ge\frac12|\partial M|$; equality in either inequality holds only for the quotient of the Nariai system, where both are attained, so $\Lambda|\partial M|<4\pi$ if $M$ contains no embedded projective plane.
\end{enumerate}
\end{thmx}

The proof of Theorem~\ref{thm:C} uses Theorem~\ref{thm:D} above. The main point is the vanishing of $E$ when $M$ is not simply connected (for simply connected $M$ it follows from Gauss's law, as recalled at the beginning of this introduction), and it rests on the fact that $H_1(M;\R)=0$ for every compact system with horizons (Theorem~\ref{thm:D}(i)). If $H_1(M;\R)\neq0$, minimizing area in a class of $H_2(M;\Z)$ not coming from $\partial M$ produces an interior locally area-minimizing surface, and Theorem~\ref{thm:B} leads to a contradiction. Hence $VE^\flat=d\psi$ with $\psi|_{\partial M}=0$, and $\int_MV|E|^2=\int_M\dv(\psi E)=0$. The area bounds in (a) and (b) then follow from the gap theorem of Ambrozio \cite[Theorem~D]{Amb17}. Once $E\equiv0$ is known, the area statement in (c) can also be derived from \cite[Theorems~B and~C]{Amb17}, applied to the universal covering; our proof passes to a finite covering with several horizons of equal area, to which Theorems~\ref{thm:D} and~\ref{thm:masscharge} apply, and it also gives the inequality for projective planes.

\begin{corollary}\label{cor:gap}
Under the assumptions of Theorem~\ref{thm:C}, either the system is the de Sitter system and $\Lambda|\partial M|=12\pi$, or $\Lambda|\partial M|\le4\pi$, with equality if and only if it is the quotient of the Nariai system.
\end{corollary}

In particular $\Lambda|\partial M|\le12\pi$, with equality only for the de Sitter system; this also follows from the Boucher--Gibbons--Horowitz--Shen inequality \cite{BGH84,Shen97,Amb17}, since $(M,g,V)$ is a static vacuum system, so the pinching hypotheses of \cite[Theorem~E]{CLS24} and \cite[Theorem~D]{FLRS25} are not needed. For several horizons, under the assumption $\sup_M|E|^2<\Lambda$, \cite[Theorem~E]{CLS24} gives an inequality weighted by the surface gravities,
\[
\sum_{j=1}^{\ell}\kappa_j\Bigl(\frac\Lambda3|\partial_jM|+\frac{16\pi^2Q_j^2}{|\partial_jM|}\Bigr)\le4\pi\sum_{j=1}^{\ell}\kappa_j,
\]
with equality only for the de Sitter system; for $\ell=1$ it reduces to $\Lambda|\partial M|\le12\pi$, which Corollary~\ref{cor:gap} improves.

\subsection{Relation to previous work}\label{sec:relation}

Theorem~\ref{thm:CLS} is a local statement, due to Ambrozio \cite[Proposition~14]{Amb17} for static vacuum systems ($E=0$), and to Cruz, Lima and Sousa \cite[Proposition~7]{CLS24} in general. Both works also used, globally in compact systems, the flow with speed $V$, which in $\inte M$ is the equidistant flow of the conformally compact metric $V^{-2}g$. In the proofs of \cite[Proposition~15]{Amb17} and \cite[Proposition~8]{CLS24}, a least area disk with boundary on $\partial M$ is flowed for all times, and $(M\setminus\partial M,g)$ is identified with the open hemisphere. In the proof of \cite[Theorem~B]{Amb17}, an interior area-minimizing surface forces the standard cylinder, that is, the Nariai system, and in Step~1 of the proof of \cite[Theorem~10]{CLS24} the flow starting at such a surface sweeps out a region diffeomorphic to $\R\times\Sph^2$, to $\RP^3$ minus a ball or to $\RP^3\#\RP^3$, which determines the topology of $M$. In these arguments the flow is continued by the compactness of minimizers in a compact manifold.

Theorems~\ref{thm:A} and~\ref{thm:B} extend these ideas. Theorem~\ref{thm:A} applies to locally area-minimizing surfaces with boundary of any topology, without assumptions on the compactness of $M$, on $\Lambda$ or on $E$; it proves that the leaves of the flow are hemispheres and identifies both the system and the surface: $M$ is the de Sitter hemisphere, and $\Sigma$ is one of its totally geodesic hemispheres, a page of the open book whose binding is the great circle $\partial\Sigma$. Theorem~\ref{thm:B} shows that $M$ is covered by a product system (Proposition~\ref{prop:globalsplit}): once a neighborhood of $\Sigma$ splits, the cross-section has curvature $K=\Lambda+E_0^2$ and the potential solves $V''=(E_0^2-\Lambda)V$ (Lemma~\ref{lem:product}), and the sign of $E_0^2-\Lambda$, the maximal positive solutions and the admissible deck groups give the models of Tables~\ref{tab:spherical} and~\ref{tab:flathyp}, eight of them with spherical cross-sections. Several of them appear as examples in \cite[Section~3]{CLS24}; for $E=0$ and $\Lambda>0$ only the Nariai system and its orientable quotient occur, and this quotient was described by Ambrozio \cite[Section~7]{Amb17}. In both theorems, analytic continuation (Section~\ref{sec:flow}) replaces the compactness of minimizers. Finally, the leaves are identified with the round hemisphere in \cite[Propositions~14(ii) and~25]{Amb17}, and with a domain bounded by geodesic arcs in a round sphere in \cite[Proposition~7(ii)]{CLS24}, through their constant curvature and a static potential vanishing on the boundary; Lemma~\ref{lem:surfhemisphere} shows directly that such a surface is a hemisphere, also when its boundary is only piecewise smooth. With a Robin-type boundary condition instead of $V^{-1}(0)=\partial M$, the same hemisphere is the equality case of a theorem of Medvedev \cite[Theorem~1.11]{Med26b}: if $R_g=6$, the boundary is minimal and $V^{-1}(0)$ is connected and meets it, then $V^{-1}(0)$ is a free boundary totally geodesic disk of area at most $2\pi$, the area of $\Sigma$ in Theorem~\ref{thm:A} for $\Lambda=3$.

For the remaining results, the homological minimization behind Theorem~\ref{thm:D}(iii) is that of the proofs of \cite[Theorem~B]{Amb17} and \cite[Theorem~10]{CLS24}; the uniqueness of the unstable horizon and the topology in Theorem~\ref{thm:D}(iv) are due to \cite[Theorem~B]{Amb17} for $E=0$ and to \cite[Theorems~A and~10]{CLS24} in general; and the area bounds in Theorem~\ref{thm:C}(a) and~(b) follow from Ambrozio's gap theorem \cite[Theorem~D]{Amb17} once $E\equiv0$ is known. Our other contributions are the vanishing of $E$ for connected horizons when $M$ is not simply connected (Theorem~\ref{thm:C}), which rests on $H_1(M;\R)=0$ (Theorem~\ref{thm:D}(i)); the area criterion for the stability of horizons, the index-one statement, the finiteness of the fundamental group and the outer-minimizing property of unions of horizons (Theorem~\ref{thm:D}(iii)--(v)); the area--charge inequality for every horizon carrying at most half of the total area, with the characterization of the charged Nariai systems as its equality case, and the bounds on the charged Hawking mass and on the charge of horizons (Theorem~\ref{thm:masscharge} and Corollary~\ref{cor:noncompactcharge}); the disk and annulus inequalities (Corollary~\ref{cor:disks}); and the results on noncompact systems (Theorem~\ref{thm:F}). For $E=0$, the finiteness of $\pi_1(M)$ is due to Ambrozio \cite[Theorem~B(i)]{Amb17}, with a different proof, and Theorem~\ref{thm:F}(iii) extends \cite[Proposition~15]{Amb17} and \cite[Proposition~8]{CLS24} to noncompact systems.

\subsection{Strategy and organization}

Section~\ref{sec:prelim} collects general properties of electrostatic systems, including their real-analyticity, Section~\ref{sec:minimization} the results on area-minimizing surfaces that we use, and Section~\ref{sec:models} the models. In Section~\ref{sec:flow} we construct the area-decreasing flow on arbitrary systems, as the projection of a Hamiltonian flow which remains smooth across the horizons, and show, by analytic continuation and a covering argument, that the whole manifold splits along a closed surface whose area is constant along the flow on one side, a hypothesis weaker than local minimality (Remark~\ref{rem:constantarea}). Section~\ref{sec:AB} proves Theorems~\ref{thm:A} and~\ref{thm:B}: for surfaces with boundary on the horizon, the problem reduces to a characterization of the three-dimensional hemisphere (Lemma~\ref{lem:hemisphere}), and for closed surfaces to the equation $V''=(E_0^2-\Lambda)V$ and to the classification of the deck groups. Sections~\ref{sec:compact} and~\ref{sec:noncompact} treat compact and noncompact systems by minimizing area in homology and homotopy classes, in $M$ and in its finite coverings, and applying Theorems~\ref{thm:A} and~\ref{thm:B} to the minimizers. Section~\ref{sec:related} contains remarks and open problems, and the Appendix a local rigidity result for area-minimizing horizons.

\section{Electrostatic systems: general properties}\label{sec:prelim}

\subsection{Conventions}
For a two-sided surface $\Sigma$ with unit normal $N$ we write $A$ for its second fundamental form, $H=\dv_\Sigma N$ for its mean curvature and $\vec H=-HN$ for its mean curvature vector; thus round spheres in $\R^3$ have mean curvature vector pointing inside. A compact domain or manifold $\Omega$ with smooth boundary is \emph{mean-convex} (resp.\ \emph{strictly mean-convex}) if $\vec H$ points into $\Omega$ or vanishes (resp.\ points into $\Omega$ and never vanishes) along $\partial\Omega$. The Jacobi operator of a two-sided minimal surface is $L_\Sigma=\Delta_\Sigma+|A|^2+\Ric(N,N)$; $\Sigma$ is \emph{stable} if $-\int_\Sigma\varphi L_\Sigma\varphi\ge0$ for all $\varphi$ (with $\varphi|_{\partial\Sigma}=0$ if $\partial\Sigma\ne\emptyset$). A number $\lambda$ is an eigenvalue of $L_\Sigma$ if $L_\Sigma\varphi+\lambda\varphi=0$ for some $\varphi\neq0$ with the same boundary condition; thus $\Sigma$ is stable if and only if the first eigenvalue of $L_\Sigma$ is non-negative, and its \emph{Morse index} is the number of negative eigenvalues of $L_\Sigma$. Along $\partial M$ we denote by $\nu$ the unit normal pointing into $M$.

\subsection{Basic identities}

\begin{lemma}\label{lem:basic}
Let $(M,g,V,E)$ be an electrostatic system with $V^{-1}(0)=\partial M$. Then:
\begin{enumerate}[label=\textup{(\roman*)}]
\item $R_g=2(\Lambda+|E|^2)$.
\item If $\partial M\neq\emptyset$, then $\nabla V\neq0$ along $\partial M$, $\nu=\nabla V/|\nabla V|$, $\partial M$ is totally geodesic, the surface gravity $\kappa=|\nabla V|$ is constant on each component of $\partial M$, and $E$ is normal to $\partial M$.
\item For every $c>0$, $(M,g,cV,E)$ is an electrostatic system.
\item The $1$-form $VE^\flat$ is closed. If $H^1(M;\R)=0$ (equivalently, $H_1(M;\R)=0$), then $VE^\flat=d\psi$ for some $\psi\in C^\infty(M)$, which is constant on each component of $\partial M$ and satisfies $V|E|^2=\langle\nabla\psi,E\rangle$ and $\dv(\psi E)=V|E|^2$.
\item If $M$ is compact and $\partial M=\partial_1M\cup\dots\cup\partial_\ell M\neq\emptyset$, with surface gravities $\kappa_j$, then
\begin{equation}\label{eq:integral}
\int_M(\Lambda-|E|^2)\,V\,dv_g=\sum_{j=1}^{\ell}\kappa_j|\partial_jM|>0 .
\end{equation}
\end{enumerate}
\end{lemma}

\begin{proof}
(i) and (ii) are \cite[Lemma~4]{CLS24}; in (ii), $\nu=\nabla V/|\nabla V|$ because $V>0$ in $\inte M$.

(iii) The system \eqref{eq:electrostatic} is linear in $V$.

(iv) If $H^1(M;\R)=0$, the closed form $VE^\flat$ is exact. Since $V=0$ on $\partial M$, the pullback of $d\psi$ to $\partial M$ vanishes, so $\psi$ is locally constant on $\partial M$. Moreover $\langle\nabla\psi,E\rangle=VE^\flat(E)=V|E|^2$ and, as $\dv E=0$, $\dv(\psi E)=\langle\nabla\psi,E\rangle$.

(v) Integrate $\Delta V=(|E|^2-\Lambda)V$ and use $\langle\nabla V,-\nu\rangle=-\kappa_j$ on $\partial_jM$.
\end{proof}

\subsection{Real-analyticity}

\begin{proposition}\label{prop:analytic}
Let $(M,g,V,E)$ be an electrostatic system. Then $\inte M$ admits a real-analytic atlas, given by harmonic coordinates, in which $g$, $V$ and $E$ are real-analytic. Consequently:
\begin{enumerate}[label=\textup{(\roman*)}]
\item every tensor field obtained from $g$, $V$, $E$ by tensor products, contractions, raising and lowering of indices and covariant differentiation is real-analytic on $\inte M$; this includes $\Riem$, $\Ric$, $R_g$, $\nabla^k\Riem$, $\Hess V$, $\nabla^kE$, $|E|^2$ and the Weyl and Cotton tensors;
\item every smooth minimal surface in $\inte M$ is real-analytic, and the normal exponential map $(s,x)\mapsto\exp_x(sN(x))$ of a two-sided minimal surface, with respect to $g$ or to $V^{-2}g$, is real-analytic where defined in $\inte M$.
\end{enumerate}
\end{proposition}

\begin{proof}
Fix $p\in\inte M$ and a simply connected neighborhood $W\subset\inte M$ of $p$ with harmonic coordinates $(x^1,x^2,x^3)$, which exist and are smooth \cite{DK81}. On $W$ we have $V>0$ and, by the Poincar\'e lemma, $VE^\flat=d\psi$ for some $\psi$, so that $E=V^{-1}\nabla\psi$. In harmonic coordinates $\Delta_g f=g^{kl}\partial_k\partial_lf$ and $\Ric_{ij}=-\frac12g^{kl}\partial_k\partial_lg_{ij}+Q_{ij}(g,\partial g)$, with $Q$ analytic in its arguments \cite{DK81}. Hence \eqref{eq:electrostatic} becomes the system
\[
\begin{aligned}
-\tfrac12g^{kl}\partial_k\partial_lg_{ij}-V^{-1}\partial_i\partial_jV&=\mathcal F_{ij}(g,\partial g,V,\partial V,\partial\psi),\\
g^{kl}\partial_k\partial_lV&=\bigl(V^{-2}g^{kl}\partial_k\psi\partial_l\psi-\Lambda\bigr)V,\\
g^{kl}\partial_k\partial_l\psi&=V^{-1}g^{kl}\partial_kV\partial_l\psi,
\end{aligned}
\]
where $\mathcal F$ is analytic for $V>0$ and the last equation is $\dv(V^{-1}\nabla\psi)=0$. Its principal symbol at $\xi\neq0$ is block triangular with diagonal blocks $-\tfrac12|\xi|_g^2\,\mathrm{Id}$, $|\xi|_g^2$ and $|\xi|_g^2$, so the system is elliptic in the sense of Douglis and Nirenberg. By Morrey's theorem \cite{Mor58}, which applies to $C^{2,\alpha}$ solutions of such systems, $(g_{ij},V,\psi)$ is real-analytic in $W$, and so is $E=V^{-1}\nabla\psi$. Harmonic functions of an analytic metric are analytic, so the transition maps between harmonic charts are analytic. This proves the first assertion (for $E=0$ and $\Lambda=0$ see \cite[Theorem~3]{MzH70}, whose remark after it also covers $\Lambda$; for $E\neq0$ and $\Lambda=0$ see \cite{Tod07}; see also \cite{CT07}). Item (i) follows because Christoffel symbols are rational functions of $g_{ij}$, $g^{ij}$ and $\partial g_{ij}$. In item (ii), the minimal surface equation is an analytic elliptic equation, so minimal surfaces are analytic \cite{Mor58}; the normal exponential map is then the flow of an analytic geodesic equation with analytic initial data, hence analytic. The metric $V^{-2}g$ is analytic on $\inte M$ since $V>0$ there.
\end{proof}

\subsection{The analytic continuation principle}

\begin{lemma}\label{lem:continuation}
Let $X$ be a connected real-analytic manifold.
\begin{enumerate}[label=\textup{(\roman*)}]
\item If a real-analytic function on $X$ vanishes on a non-empty open set, it vanishes identically.
\item If two real-analytic tensor fields on $X$ coincide on a non-empty open set, they coincide on $X$.
\item If $F:X\to\inte M$ is real-analytic and $T$ is a real-analytic tensor field on $\inte M$, then $F^*T$ is real-analytic; in particular, if $F^*g$, $V\circ F$ or $F^*E$ coincide on a non-empty open set with analytic objects defined on $X$, they coincide on $X$.
\end{enumerate}
\end{lemma}

\begin{proof}
(i) This is the identity theorem; we recall the argument. Let $f$ be such a function and $Z$ the set of points of $X$ having a neighborhood on which $f$ vanishes. It is open by definition and non-empty by hypothesis. Let $p$ be in the closure of $Z$ and fix an analytic chart around $p$. All partial derivatives of $f$ vanish on $Z$, hence at $p$ by continuity, so the Taylor series of $f$ at $p$ is zero. Since $f$ is real-analytic, it is the sum of this series in a neighborhood of $p$; thus $f$ vanishes there, and $p\in Z$. So $Z$ is closed, and $Z=X$ because $X$ is connected.

(ii) Let $Z$ be the set of points having a neighborhood on which the two fields coincide. In an analytic chart, the components of their difference are real-analytic functions, so the argument of (i) shows that $Z$ is open, closed and non-empty.

(iii) In analytic charts, the components of $F^*T$ are polynomials in the first derivatives of $F$ whose coefficients are the components of $T$ composed with $F$, hence real-analytic. The last assertion follows from (i) and (ii).
\end{proof}

\subsection{Product systems}

\begin{definition}\label{def:product}
Let $(\Sigma,h)$ be a closed orientable surface, $I\subset\R$ an interval, $E_0\in\R$ and $V\in C^\infty(I)$. The quadruple $(I\times\Sigma,\,ds^2+h,\,V(s),\,E_0\partial_s)$ is a \emph{product system} if $V>0$ in the interior of $I$, $V$ vanishes exactly at the endpoints belonging to $I$, and \eqref{eq:electrostatic} holds. An electrostatic system is \emph{of product type} if it is isometric to $(I\times\Sigma)/\Gamma$, where $(I\times\Sigma,ds^2+h,V,E_0\partial_s)$ is a product system and $\Gamma$ is a group of isometries acting freely and properly discontinuously and preserving $V$ and $E_0\partial_s$.
\end{definition}

\begin{lemma}\label{lem:product}
Let $(\Sigma,h)$ be a closed surface, $I$ an interval, $V=V(s)>0$ and $E$ a vector field on $I\times\Sigma$ such that $(I\times\Sigma,ds^2+h,V,E)$ satisfies \eqref{eq:electrostatic}. Then $E=E_0\partial_s$ for a constant $E_0$, $h$ has constant Gaussian curvature $K=\Lambda+E_0^2$, and $V''=(E_0^2-\Lambda)V$. Conversely, these conditions imply \eqref{eq:electrostatic}.
\end{lemma}

\begin{proof}
For the product metric $g=ds^2+h$ we have $\Ric_g=K_h\,h$ and $\Hess V=V''ds^2$. Write $E=f\partial_s+E^T$ with $E^T$ tangent to $\Sigma$. The tangential part of the first equation in \eqref{eq:electrostatic} reads $K_hh=(\Lambda+|E|^2)h-2(E^T)^\flat\otimes(E^T)^\flat$; since $h$ has rank two, $E^T=0$ and $K_h=\Lambda+f^2$. Then $\dv E=\partial_sf=0$ and $d(VE^\flat)=V\,d_\Sigma f\wedge ds=0$, so $f=E_0$ is constant, $K_h=\Lambda+E_0^2$ is constant, and the $ss$-component reads $0=V''+(\Lambda-E_0^2)V$. The mixed components vanish identically. The converse is the same computation.
\end{proof}

\begin{lemma}\label{lem:noconvex}
Let $(M,g)$ be a Riemannian $3$-manifold containing an isometrically embedded product region $P=[a,b]\times\Sigma$, $g|_P=ds^2+h$, with $\Sigma$ closed and connected. If $\Omega\subset M$ is a compact domain with smooth boundary such that $\{a\}\times\Sigma\subset\inte\Omega$ and $(\{b\}\times\Sigma)\cap\Omega=\emptyset$, then at some point of $\partial\Omega\cap P$ the mean curvature vector of $\partial\Omega$ points out of $\Omega$ or vanishes; in particular, $\Omega$ is not strictly mean-convex. Hence a noncompact manifold having an end isometric to $((a,\infty)\times\Sigma,ds^2+h)$ does not have strictly mean-convex ends.
\end{lemma}

\begin{proof}
Let $t_0=\sup\{t\in[a,b]:[a,t]\times\Sigma\subset\Omega\}<b$. The open set $(a,t_0)\times\Sigma$ is contained in $\Omega$, hence in $\inte\Omega$, and by compactness of $\Sigma$ there is $q\in(\{t_0\}\times\Sigma)\cap\partial\Omega$. Near $q$ the surface $\partial\Omega$ lies in $\{s\ge t_0\}$ and touches the totally geodesic slice $\{t_0\}\times\Sigma$ at $q$, so locally $\partial\Omega=\{s=t_0+u(x)\}$ with $u\ge0$, $u(q)=0$, $\nabla u(q)=0$ and $\Hess_hu(q)\ge0$, while $\Omega$ lies below the graph. Thus near $q$ the surface $\partial\Omega$ is the level set $\{F=t_0\}$ of $F(s,x)=s-u(x)$, and $\nabla F(q)=\partial_s$. Since $g=ds^2+h$, we have $\Hess_gs=0$ and $\Hess_g(u\circ\pi)=\pi^*\Hess_hu$, where $\pi(s,x)=x$, so the second fundamental form of $\partial\Omega$ at $q$ with respect to $\partial_s=\nabla F/|\nabla F|$ is the restriction of $\Hess_gF/|\nabla F|$ to $T_q\partial\Omega=T_q(\{t_0\}\times\Sigma)$, that is, $-\Hess_hu(q)$, so $\vec H(q)=(\Delta_hu)(q)\,\partial_s$ points out of $\Omega$ or vanishes. For the last assertion, given an exhaustion $\{\Omega_k\}$, choose $k$ with $\{a+1\}\times\Sigma\subset\inte\Omega_k$ and $b>a+1$ with $(\{b\}\times\Sigma)\cap\Omega_k=\emptyset$; since the end lies in $\inte M$, the point given by the first part lies on $\partial_+\Omega_k$.
\end{proof}

\section{Area-minimizing surfaces}\label{sec:minimization}

In this section $W$ denotes a compact connected orientable Riemannian $3$-manifold, possibly with boundary. We collect the existence and regularity results used in the paper.

\subsection{Definitions}

\begin{definition}[{\cite[Definition~6]{CLS24}}]\label{def:almost}
A compact immersed surface $\Sigma\subset M$ is \emph{almost properly embedded} if $\Sigma\cap\partial M=\partial\Sigma$ and there are disjoint finite sets $\{p_i\},\{q_j\}\subset\partial\Sigma$ such that $\Sigma\setminus\{p_i\}$ is embedded and $\Sigma\setminus\{q_j\}$ is smooth. The points near which $\Sigma$ is not a smooth embedded surface form its \emph{exceptional set}, and $\Sigma$ is \emph{properly embedded} if this set is empty.
\end{definition}

We assume in addition that almost properly embedded surfaces have finite area, as the least area disks and annuli of Theorems~\ref{thm:plateau} and~\ref{thm:douglas} do; this is used in Lemma~\ref{lem:totgeod} and to compare areas along the flow in Section~\ref{sec:flow}.

\begin{definition}\label{def:lam}
Let $\Sigma\subset M$ be a compact two-sided almost properly embedded surface with exceptional set $S$; for instance, $\Sigma$ may be a closed embedded surface in $\inte M$, for which $S=\emptyset$. A \emph{deformation} of $\Sigma$ is a smooth map $F:[0,1]\times(\Sigma\setminus S)\to M$ such that $F(0,\cdot)$ is the inclusion and each $F(s,\cdot)$ is an immersion which is the identity on $\partial\Sigma\setminus S$ and extends continuously to $\Sigma$ by the identity on $S$. Its \emph{size} is $\sup_s\|F(s,\cdot)-F(0,\cdot)\|_{C^1}$, measured in a fixed finite atlas of a neighborhood of $\Sigma$.
\begin{enumerate}[label=\textup{(\roman*)}]
\item $\Sigma$ is \emph{locally area-minimizing} if there is $\varepsilon>0$ such that $|\Sigma|\le|F(1,\Sigma\setminus S)|$ for every deformation $F$ of $\Sigma$ of size less than $\varepsilon$, the area of the immersed surface $F(1,\Sigma\setminus S)$ being counted with multiplicity.
\item A one-sided closed surface $\Sigma\subset\inte M$ is locally area-minimizing if the two-sided preimage $\hat\Sigma$ of $\Sigma$ in the double cover $\hat M\to M$ associated with the mod $2$ intersection number with $\Sigma$ is locally area-minimizing in $\hat M$.
\end{enumerate}
\end{definition}

Every locally area-minimizing surface is minimal and stable for variations with compact support in its interior: for a unit normal $N$ and $\varphi\in C^\infty_c(\inte\Sigma)$, $F(s,x)=\exp_x(s\varphi(x)N(x))$ is a deformation of size at most $C\|\varphi\|_{C^1}$ once $\varphi$ is small. If $S=\emptyset$, deformations of small size are isotopies through embeddings; immersions are allowed because the leaves of the flow of Section~\ref{sec:flow} may be only immersed near the exceptional set (Lemma~\ref{lem:flow}). For closed surfaces, (i) is weaker than the notion of local least area of Cai and Galloway \cite[Section~2]{CG00}, in which $\Sigma$ is compared with all surfaces isotopic to it in a normal neighborhood; compare also the minimization in isotopy classes of \cite{MSY82} and, for free boundary surfaces, the comparison with all nearby properly immersed surfaces in \cite{Amb15}. In~(ii), $\hat\Sigma$ is the orientable double cover of $\Sigma$ and $(\hat M,\hat g,\hat V,\hat E)$, with the pulled back data, is an electrostatic system.

\subsection{The maximum principle}

\begin{lemma}\label{lem:maxprinciple}
Let $\Omega$ be a compact mean-convex domain and let $\Sigma\subset\Omega$ be a connected minimal surface, possibly with boundary. If an interior point of $\Sigma$ lies on $\partial\Omega$, then $\Sigma$ is contained in a component of $\partial\Omega$ which is minimal near $\Sigma$; in particular, $\inte\Sigma$ does not meet strictly mean-convex parts of $\partial\Omega$. If $\Sigma$ is closed and meets a totally geodesic component $B$ of $\partial\Omega$, then $\Sigma=B$.
\end{lemma}

\begin{proof}
At a contact point $\Sigma$ and $\partial\Omega$ are tangent and, locally, $\Sigma$ lies on the mean-convex side of $\partial\Omega$; writing both as graphs over the common tangent plane, the strong maximum principle for the mean curvature operator \cite{SW89,Whi10} shows that they coincide near the contact point, which is impossible where $\vec H_{\partial\Omega}\neq0$. The set of interior points of $\Sigma$ lying on $\partial\Omega$ is therefore open and closed in $\inte\Sigma$. If $\Sigma$ is closed and $\Sigma\subset B$, then $\Sigma=B$ since both are closed connected surfaces.
\end{proof}

\subsection{Minimization in homology classes}

\begin{theorem}\label{thm:homology}
Let $W$ be mean-convex and let $\alpha\in H_2(W;\Z)$. Then there are pairwise disjoint, closed, connected, embedded, orientable surfaces $\Sigma_1,\dots,\Sigma_R\subset W$ and positive integers $n_1,\dots,n_R$ with $\sum_in_i[\Sigma_i]=\alpha$ and
\[
\sum_in_i|\Sigma_i|=\inf\{\mathbf M(T):T\text{ an integral $2$-cycle in }W\text{ representing }\alpha\}.
\]
Each $\Sigma_i$ either is a minimal component of $\partial W$ or is contained in $\inte W$, and in the latter case it is locally area-minimizing. No $\Sigma_i$ meets a strictly mean-convex component of $\partial W$.
\end{theorem}

\begin{proof}
Existence of a mass minimizing integral current $T$ in $\alpha$ follows from the Federer--Fleming compactness theorem \cite[4.4.4]{Fed69}; the mean-convex boundary acts as a barrier, so $T$ may be found by minimizing in a slightly larger manifold and the maximum principle confines its support to $W$ \cite[Section~6]{MSY82}, \cite{Whi10}. In dimension three, codimension one area-minimizing currents are smooth embedded surfaces away from $\partial W$ \cite[5.4.15]{Fed69}. Near $\partial W$, the maximum principle for stationary varifolds \cite{Whi10} shows that the support of $T$ does not meet the strictly mean-convex part of $\partial W$, and the strong maximum principle \cite{SW89} shows that the support contains every minimal component of $\partial W$ that it touches. Hence each component of the support either lies in $\inte W$ or is a minimal component of $\partial W$ (for smooth surfaces this is also Lemma~\ref{lem:maxprinciple}). Thus $T=\sum n_i[\Sigma_i]$, where $\Sigma_1,\dots,\Sigma_R$ are the components of the support of $T$, which are closed, embedded and oriented, and $n_i$ are positive integers \cite[4.1.31 and 5.4.15]{Fed69}; see also \cite[proof of Theorem~B(ii)]{Amb17} and \cite[proof of Theorem~10]{CLS24}. Each $\Sigma_i$ is two-sided because $W$ is orientable. If $\Sigma_i\subset\inte W$ and $F$ is a deformation of $\Sigma_i$ in $W$ (Definition~\ref{def:lam}), the $3$-chain $F_\#([0,1]\times[\Sigma_i])$ has boundary $F(1,\cdot)_\#[\Sigma_i]-[\Sigma_i]$, so replacing $n_i[\Sigma_i]$ by $n_iF(1,\cdot)_\#[\Sigma_i]$ yields a competitor; hence $|\Sigma_i|\le|F(1,\Sigma_i)|$.
\end{proof}

\begin{lemma}\label{lem:boundaryclasses}
Let $W$ be compact, connected and oriented, and let $B_1,\dots,B_m$ be the components of $\partial W$, with the boundary orientation. Then the kernel of $\bigoplus_j\Z[B_j]\to H_2(W;\Z)$ is generated by $\sum_j[B_j]$. In particular, if $m\ge2$, the classes of any proper subfamily of $\{B_j\}$ are linearly independent in $H_2(W;\R)$, and $[B_1]\notin\operatorname{span}\{[B_j]:j\in J\}$ whenever $1\notin J$ and $J\cup\{1\}\neq\{1,\dots,m\}$.
\end{lemma}

\begin{proof}
This follows from the exact sequence $H_3(W,\partial W;\Z)\to H_2(\partial W;\Z)\to H_2(W;\Z)$, since $H_3(W,\partial W;\Z)\cong\Z$ is generated by the fundamental class, whose boundary is $\sum_j[B_j]$.
\end{proof}

\subsection{Disks, loops and annuli}

\begin{theorem}[\cite{MY80,MY82a,MY82b}]\label{thm:plateau}
Let $W$ be mean-convex and let $\gamma\subset\partial W$ be a smooth Jordan curve which is null-homotopic in $W$. Then $\gamma$ bounds a disk $D$ of least area among all disks in $W$ bounded by $\gamma$, and every such disk is almost properly embedded or contained in $\partial W$.
\end{theorem}

The alternative $D\subset\partial W$ can only occur when $\gamma$ bounds a disk in a minimal component of $\partial W$; this follows from Lemma~\ref{lem:maxprinciple}. For mean-convex boundary, the existence and the embeddedness of $D$ are \cite[Theorem~1]{MY82b}, and, by \cite[Theorem~2]{MY82b}, a least area disk which is not contained in $\partial W$ meets $\partial W$ only along $\gamma$ and has no boundary branch points, so it is in fact properly embedded. Lemma~\ref{lem:flow} and Theorem~\ref{thm:A} are nevertheless stated for the larger class of almost properly embedded surfaces of \cite[Definition~6]{CLS24}.

\begin{theorem}[Loop theorem \cite{Pap57,Hem76}]\label{thm:loop}
Let $W$ be a compact $3$-manifold and $F$ a component of $\partial W$ such that $\ker(\pi_1(F)\to\pi_1(W))\neq1$. Then there is a properly embedded disk $D\subset W$ whose boundary is a Jordan curve in $F$ which is not null-homotopic in $F$.
\end{theorem}

\begin{theorem}[\cite{Dou31,Jos85,MY82a,MY82b}]\label{thm:douglas}
Let $W$ be mean-convex and let $\Gamma_1,\Gamma_2\subset\partial W$ be disjoint smooth Jordan curves, null-homotopic in $W$. Every annulus of least area among the annuli in $W$ bounded by $\Gamma_1\cup\Gamma_2$ is either almost properly embedded or contained in $\partial W$, and such an annulus exists if the strict Douglas condition holds in $W$.
\end{theorem}

For convex boundary, the existence is \cite[Theorem~1]{MY82a}, and every least area annulus is embedded by \cite[Theorem~5]{MY82a}: joining two disks bounded by $\Gamma_1$ and $\Gamma_2$ by a thin tube shows that the infimum of the areas of the annuli is at most that of the pairs of disks, so a least area annulus also has least area among all planar domains, possibly disconnected, bounded by $\Gamma_1\cup\Gamma_2$, which is the class considered there. The arguments of \cite[Section~1]{MY82b} extend these results to mean-convex boundary, and, by \cite[Theorem~2]{MY82b}, a least area annulus which is not contained in $\partial W$ meets $\partial W$ only along $\Gamma_1\cup\Gamma_2$ and has no boundary branch points. The curves are assumed null-homotopic so that the disks in the strict Douglas condition exist.

\subsection{Nodal sets}

\begin{lemma}\label{lem:nodal}
Let $(S,h)$ be a Riemannian surface diffeomorphic to $\Sph^2$, $q\in C^\infty(S)$, and $u$ an eigenfunction of $\Delta_h+q$ with exactly two nodal domains. Then $0$ is a regular value of $u$, the nodal set $\Gamma=u^{-1}(0)$ is a smooth Jordan curve, and the nodal domains are disks $\Omega_1,\Omega_2$ with $\partial\Omega_i=\Gamma$, on which $|u|$ is a first Dirichlet eigenfunction.
\end{lemma}

\begin{proof}
By Cheng's structure theorem \cite{Che76}, $\Gamma$ is a finite graph whose vertices are the critical points of $u$ on $\Gamma$, each vertex having even degree at least four, and $\Gamma$ has no endpoints. If $\Gamma$ had a vertex, or two components, then Euler's formula on $\Sph^2$ would give at least three faces, i.e.\ nodal domains. Hence $\Gamma$ is a single smooth closed curve, and the two complementary components are disks by the Jordan--Schoenflies theorem. The last assertion is standard \cite{CH53}.
\end{proof}

\section{Standard models}\label{sec:models}

The \emph{de Sitter system} is the closed hemisphere $\Sph^3_+(K)$ of constant curvature $K=\Lambda/3$ with $V$ the height function and $E=0$; its boundary is a totally geodesic sphere of Morse index one. The other models are of product type.

By Lemma~\ref{lem:product}, a product system $(I\times\Sigma,ds^2+h,V(s),E_0\partial_s)$ is determined by $\Lambda$, $E_0$, a closed surface $(\Sigma,h)$ of constant curvature $K=\Lambda+E_0^2$, and a positive solution of $V''=\mu V$, $\mu=E_0^2-\Lambda$, on the interior of $I$, vanishing at the endpoints of $I$ that belong to $I$; these endpoints are the horizons, and the surface gravity there is $\kappa=|V'|$. Up to a translation or reflection of $s$, the maximal positive solutions are the positive multiples of
\[
\begin{array}{ll}
\mu>0: & \cosh(\sqrt\mu\,s)\text{ on }\R,\quad \sinh(\sqrt\mu\,s)\text{ on }[0,\infty),\quad e^{\sqrt\mu\,s}\text{ on }\R;\\[2pt]
\mu=0: & 1\text{ on }\R,\quad s\text{ on }[0,\infty);\\[2pt]
\mu<0: & \cos(\sqrt{-\mu}\,s)\text{ on }[-a,a],\quad a=\pi/(2\sqrt{-\mu}).
\end{array}
\]
A system of product type is a quotient of a product system by a group of orientation-preserving isometries of the form $(s,x)\mapsto(\sigma(s),\varphi(x))$ preserving $V$ and $E_0\partial_s$ (see the proof of Theorem~\ref{thm:B}); in particular, $\sigma$ is a translation unless $E_0=0$.

\subsection{Spherical cross-sections}

When $K>0$, the cross-section is the round sphere of curvature $K$. Since every rotation of $\Sph^2$ has a fixed point, a deck transformation is determined by its action on $s$, and the proof of Theorem~\ref{thm:B} gives the eight systems of Table~\ref{tab:spherical}. Here $\mu>0$ means $E_0^2>|\Lambda|$, and $\mu\le0$ forces $\Lambda>0$.

\begin{table}[ht]
\centering\small
\renewcommand{\arraystretch}{1.25}
\begin{tabular}{llllll}
\hline
 & System & $\mu$ & $V(s)$ & $M$ & Horizons ($\kappa$) \\
\hline
(i) & global Bertotti--Robinson & $>0$ & $\cosh(\sqrt\mu\,s)$ & $\R\times\Sph^2$ & none \\
(ii) & cold black hole & $>0$ & $\sinh(\sqrt\mu\,s)$ & $[0,\infty)\times\Sph^2$ & one ($\sqrt\mu$) \\
(iii) & Poincar\'e Bertotti--Robinson & $>0$ & $e^{\sqrt\mu\,s}$ & $\R\times\Sph^2$ & none \\
(iv) & Nariai--Bertotti--Robinson & $0$ & $1$ & $\R\times\Sph^2$ & none \\
(v) & compact quotients of (iv) & $0$ & $1$ & $\Sph^1\times\Sph^2$ & none \\
(vi) & ultracold black hole & $0$ & $s$ & $[0,\infty)\times\Sph^2$ & one ($1$) \\
(vii) & charged Nariai & $<0$ & $\cos(\sqrt{-\mu}\,s)$ & $[-a,a]\times\Sph^2$ & two ($\sqrt{-\mu}$) \\
(viii) & quotient of Nariai & $-\Lambda$ & $\cos(\sqrt\Lambda\,s)$ & $\RP^3\setminus B^3$ & one ($\sqrt\Lambda$) \\
\hline
\end{tabular}
\caption{Systems of product type with $K=\Lambda+E_0^2>0$.}\label{tab:spherical}
\end{table}

The system (iv) is also called the Pleba\'nski--Hacyan system \cite{PH79,GP09}, and (v) is its quotient by $(s,p)\mapsto(s+L,Ap)$, with $L>0$ and $A\in SO(3)$; since $SO(3)$ is connected, $M\cong\Sph^1\times\Sph^2$, and for $A\neq\mathrm{Id}$ the quotient is not isometric to the Riemannian product $\Sph^1_L\times\Sph^2$. The charged Nariai systems (vii) reduce to the Nariai system for $E_0=0$, and (viii) is the quotient of the Nariai system by the free orientation-preserving involution $(s,p)\mapsto(-s,-p)$, described by Ambrozio \cite[Section~7]{Amb17}; the image of $\{0\}\times\Sph^2$ is a one-sided projective plane, and $\partial M\cong\Sph^2$. All these horizons are stable. Systems (i)--(iii) exist for every sign of $\Lambda$, and for $\Lambda=0$, (i) is the Bertotti--Robinson system proper.

\subsection{Flat and hyperbolic cross-sections}

When $K\le0$, the cross-section is a closed orientable surface $\Sigma$ of positive genus with constant curvature $K$: a flat torus $T^2$ if $K=0$, and a hyperbolic surface $\Sigma_{\mathrm g}$ of genus $\mathrm g\ge2$ if $K<0$. Then $\Lambda=K-E_0^2\le0$ and $\mu=2E_0^2-K\ge0$, so the charged Nariai systems have no analogue. Isometries of $\Sigma$ may now act freely, so deck transformations that fix $s$ need not be trivial; absorbing them into $\Sigma$, the argument in the proof of Theorem~\ref{thm:B} gives the systems of Table~\ref{tab:flathyp}, labelled after their spherical analogues. Thus (vii$'$) does not occur, and neither does (v$'$), since the compact quotients of $\R\times T^2$ are included in (iv$'$). Like (viii), (viii$'$) is a quotient by a reflection of $s$ composed with a free orientation-reversing involution of the cross-section, which forces $E_0=0$; the covering system is then (i$'$), since $\mu=-K>0$.

\begin{table}[ht]
\centering\small
\renewcommand{\arraystretch}{1.25}
\begin{tabular}{llllll}
\hline
 & System & $\mu$ & $V(s)$ & $M$ & Horizons ($\kappa$) \\
\hline
(i$'$) & global $AdS_2\times\Sigma$ & $>0$ & $\cosh(\sqrt\mu\,s)$ & $\R\times\Sigma$ & none \\
(ii$'$) & $AdS_2$ black hole $\times\,\Sigma$ & $>0$ & $\sinh(\sqrt\mu\,s)$ & $[0,\infty)\times\Sigma$ & one ($\sqrt\mu$) \\
(iii$'$) & Poincar\'e $AdS_2\times\Sigma$ & $>0$ & $e^{\sqrt\mu\,s}$ & $\R\times\Sigma$ & none \\
(iv$'$) & flat systems & $0$ & $1$ & $(\R\times T^2)/\Gamma$ & none \\
(vi$'$) & Rindler-type & $0$ & $s$ & $[0,\infty)\times T^2$ & one ($1$) \\
(viii$'$) & quotient of (i$'$), $E_0=0$ & $-\Lambda$ & $\cosh(\sqrt{-\Lambda}\,s)$ & $(\R\times\Sigma_{\mathrm g})/\langle\tau\rangle$ & none \\
\hline
\end{tabular}
\caption{Systems of product type with $K=\Lambda+E_0^2\le0$.}\label{tab:flathyp}
\end{table}

Here $\Sigma$ is $T^2$ or $\Sigma_{\mathrm g}$. The systems (i$'$)--(iii$'$) come from static regions of $AdS_2\times\Sigma$: with hyperbolic cross-sections they are the anti-Nariai systems, charged if $E_0\neq0$, and with flat cross-sections, where $E_0\neq0$ and $\Lambda=-E_0^2<0$, the Pleba\'nski--Hacyan-type systems \cite{GP09,CM25}. In (iv$'$), $E=0$ and $\Lambda=0$, and $\Gamma$ is a discrete group of orientation-preserving isometries of $\R\times T^2$ of the above form; these are flat systems, and the compact ones are the closed flat manifolds of Theorem~\ref{thm:E}(a). The Rindler-type system (vi$'$) also has $E=0$ and $\Lambda=0$. In (viii$'$), $\tau(s,x)=(-s,\iota(x))$ for a free orientation-reversing isometric involution $\iota$ of $\Sigma_{\mathrm g}$, which exists when $\Sigma_{\mathrm g}$ double covers a non-orientable surface isometrically. In particular, no system of product type with $K\le0$ is compact with non-empty boundary.

\begin{remark}[Systems of product type]\label{rem:sharp}
A noncompact system of product type has an end isometric to $((s_0,\infty)\times\Sigma',ds^2+h')$ with $\Sigma'$ closed, so, by Lemma~\ref{lem:noconvex}, it does not have strictly mean-convex ends. Its ends are mean-convex: the systems (ii), (vi), (ii$'$) and (vi$'$) are exhausted by $\Omega_k=[0,k]\times\Sigma$, and the cylinders (i), (iii), (iv), (i$'$) and (iii$'$) by $[-k,k]\times\Sigma$, with totally geodesic outer boundaries. These systems show that strict mean-convexity is needed in Theorem~\ref{thm:F}(i) and~(ii): in the first four, $\partial M$ is stable, since its Jacobi operator is $\Delta$ because $\Ric(\partial_s,\partial_s)=0$, but not strictly outer-minimizing, and the cylinders have two ends and $H_2(M;\Z)\cong\Z$. On the other hand, the proofs in Section~\ref{sec:noncompact} use the strict mean-convexity of the ends only to exclude systems of product type. Thus Theorem~\ref{thm:F}(i) and~(ii) and Corollary~\ref{cor:noncompactcharge} hold for every noncompact system with mean-convex ends which is not of product type, and Theorem~\ref{thm:F}(iii) and~(iv) for every noncompact system with mean-convex ends.
\end{remark}

\section{The area-decreasing flow and the global splitting}\label{sec:flow}

\subsection{Short-time existence}

\begin{lemma}\label{lem:flow}
Let $(M,g,V,E)$ be an electrostatic system with $V^{-1}(0)=\partial M$ and $\Sigma_0\subset M$ a compact, connected, two-sided surface with unit normal $N_0$, which is either closed and contained in $\inte M$, or almost properly embedded with $\Sigma_0\cap\partial M=\partial\Sigma_0\neq\emptyset$, exceptional set $S$, and bounded second fundamental form on $\Sigma_0\setminus S$; we set $S=\emptyset$ in the first case. Then there are $\delta>0$ and a smooth map $\Phi:[0,\delta)\times(\Sigma_0\setminus S)\to M$ such that $\Phi(0,\cdot)$ is the inclusion, each $\Phi(t,\cdot)$ is an immersion, and an embedding if $S=\emptyset$, which fixes $\partial\Sigma_0\setminus S$, maps $\inte\Sigma_0$ into $\inte M$ and extends continuously to $\Sigma_0$ by $\Phi(t,q)=q$ for $q\in S$, and
\[
\partial_t\Phi=(V\circ\Phi)\,N_\Phi ,
\]
where $N_\Phi$ is the unit normal of $\Sigma_t=\Phi(t,\Sigma_0)$ extending $N_0$. The flow depends only on $(g,V)$ near $\Sigma_0$; consequently, Theorem~\ref{thm:CLS} holds without assuming that $M$ is compact.
\end{lemma}

\begin{proof}
If $\partial\Sigma_0\neq\emptyset$, we extend $g$ and $V$ smoothly across $\partial M$ near $\partial\Sigma_0$; since $\nabla V\neq0$ on $\partial M$ (Lemma~\ref{lem:basic}(ii)), the extension of $V$ is negative outside $M$ near $\partial M$. On the cotangent bundle of $M$, extended in this way, minus its zero section, consider the Hamiltonian
\[
\mathcal H(q,p)=V(q)\,|p|_{g(q)} ,
\]
which is smooth, also where $V=0$, and homogeneous of degree one in $p$. Let $(q,p)(t,x)$ be the solution of Hamilton's equations
\[
\dot q=\partial_p\mathcal H=V(q)\,\frac{p^\sharp}{|p|},\qquad\dot p=-\partial_q\mathcal H,\qquad q(0,x)=x,\quad p(0,x)=N_0(x)^\flat ,
\]
for $x\in\Sigma_0$. Since $\Sigma_0$ is compact and smooth up to its boundary, the solutions exist on a common interval $[0,\delta_1)$ and depend smoothly on $(t,x)$; set $\Phi(t,x)=q(t,x)$.

Let $X$ be a vector field on $\Sigma_0$. Differentiating Euler's identity $\langle p,\partial_p\mathcal H\rangle=\mathcal H$ in $q$ and in $p$ gives $\langle p,\partial_X(\partial_p\mathcal H)\rangle=\langle\partial_q\mathcal H,\partial_Xq\rangle$, hence
\[
\partial_t\langle p,\partial_Xq\rangle=-\langle\partial_q\mathcal H,\partial_Xq\rangle+\langle p,\partial_X(\partial_p\mathcal H)\rangle=0 .
\]
Thus $p(t,x)$ annihilates $d\Phi(t,\cdot)(T_x\Sigma_0)$, as it does at $t=0$. For small $t$, $\Phi(t,\cdot)$ is $C^1$-close to the inclusion, hence an embedding onto a surface $\Sigma_t$, whose unit normal extending $N_0$ is $N_\Phi=p^\sharp/|p|$; therefore $\partial_t\Phi=(V\circ\Phi)N_\Phi$. Finally, $\partial_t(V\circ\Phi)=(V\circ\Phi)\langle\nabla V,N_\Phi\rangle$, so $V\circ\Phi$ vanishes identically or has no zeros along each curve $t\mapsto\Phi(t,x)$. Hence $\Phi(t,x)=x$ for $x\in\partial\Sigma_0$, since $\dot q=0$ where $V=0$, while $V(\Phi(t,x))>0$ for $x\in\inte\Sigma_0$, so that $\Phi(t,\inte\Sigma_0)\subset\inte M$ for small $t$. In particular $\Sigma_t\cap\partial M=\partial\Sigma_0$. In $\inte M$, $\mathcal H$ is the norm of the cometric of $V^{-2}g$, so $\Phi(t,x)$ is the normal exponential map of $V^{-2}g$ at time $t$, as in \cite[Proposition~7]{CLS24}; the point of the Hamiltonian formulation is that it remains smooth across $V^{-1}(0)$.

If $\Sigma_0$ is almost properly embedded, with exceptional set $S$, the argument applies on $\Sigma_0\setminus S$; only two steps used the compactness of $\Sigma_0$ and its smoothness up to the boundary. First, since $\mathcal H$ is homogeneous of degree one in $p$, the solution with initial data $(x,\lambda\xi)$, $\lambda>0$, is $(q,\lambda p)$, where $(q,p)$ is the solution with initial data $(x,\xi)$. As the unit covectors over the compact set $\Sigma_0$ form a compact set, the solutions with initial data $(x,N_0(x)^\flat)$, $x\in\Sigma_0\setminus S$, exist on a common interval $[0,\delta_1)$, and, in coordinates from a finite atlas around $\Sigma_0$, the time-$t$ map $\Theta_t$ of the Hamiltonian flow satisfies $d\Theta_t=I+O(t)$ uniformly on these unit covectors. Second, $|\nabla N_0|=|A|$ is bounded, so $\iota(x)=(x,N_0(x)^\flat)$ has bounded derivative. Since $\Phi(t,\cdot)=\pi\circ\Theta_t\circ\iota$, with $\pi(q,p)=q$, and $d\pi\circ d\iota$ is the identity, we get, in these coordinates, $|d\Phi(t,\cdot)_x(v)-v|\le C_1t|v|$ for $x\in\Sigma_0\setminus S$ and $v\in T_x\Sigma_0$, and therefore
\[
|d\Phi(t,\cdot)_x(v)|\ge|v|-|d\Phi(t,\cdot)_x(v)-v|\ge(1-C_1t)\,|v| .
\]
Hence $\Phi(t,\cdot)$ is an immersion of $\Sigma_0\setminus S$ for $t<\min\{\delta_1,1/C_1\}$, with unit normal $N_\Phi=p^\sharp/|p|$.

To extend $\Phi$ to $S$, we take $\delta\le\min\{\delta_1,1/C_1\}$ so small that the curves $\tau\mapsto\Phi(\tau,x)$, $\tau\in[0,\delta)$, $x\in\Sigma_0\setminus S$, stay in a compact neighborhood of $\Sigma_0$, and we let $C_2$ bound $|\nabla V|$ on it. Since $V\circ\Phi\ge0$ and $\partial_t(V\circ\Phi)=(V\circ\Phi)\langle\nabla V,N_\Phi\rangle\le C_2\,V\circ\Phi$, we have $V(\Phi(\tau,x))\le e^{C_2\tau}V(x)$. As $|\partial_\tau\Phi|=V\circ\Phi$, the curve $\tau\mapsto\Phi(\tau,x)$, $\tau\in[0,t]$, gives
\begin{equation}\label{eq:distflow}
\begin{aligned}
\operatorname{dist}_g\bigl(\Phi(t,x),x\bigr)&\le\int_0^t|\partial_\tau\Phi(\tau,x)|\,d\tau=\int_0^tV(\Phi(\tau,x))\,d\tau\\
&\le\int_0^te^{C_2\tau}V(x)\,d\tau\le te^{C_2t}\,V(x)
\end{aligned}
\end{equation}
for $x\in\Sigma_0\setminus S$ and $t\in[0,\delta)$. If $q\in S$, then $V(q)=0$, because $S\subset\partial\Sigma_0\subset\partial M$, and \eqref{eq:distflow} gives
\[
\operatorname{dist}_g\bigl(\Phi(t,x),q\bigr)\le\operatorname{dist}_g(x,q)+\delta e^{C_2\delta}\,V(x)\longrightarrow0\qquad\text{as }x\to q,
\]
uniformly in $t$. Hence $\Phi$ extends continuously to $[0,\delta)\times\Sigma_0$ by $\Phi(t,q)=q$ for $q\in S$.

It remains to check that Theorem~\ref{thm:CLS} holds without assuming that $M$ is compact; its proof only involves $(g,V,E)$ on $U=\Phi([0,\delta)\times\Sigma_0)$. With our convention $H=\dv_\Sigma N$, the variation formula for the mean curvature under a normal variation with speed $V$, together with \eqref{eq:electrostatic}, gives
\[
\partial_tH_t=-L_{\Sigma_t}V=\langle\nabla V,N_t\rangle H_t-\bigl(|A_t|^2+2|E^T|^2\bigr)V
\]
on $\Sigma_t\setminus S$, where $E^T$ is the component of $E$ tangent to $\Sigma_t$; compare \cite[proof of Proposition~7]{CLS24}. The rest of that proof uses only this equation and \eqref{eq:electrostatic} on $U$, applied on $\Sigma_0\setminus S$ with $S$ finite, so the compactness of $M$ is not needed.
\end{proof}

\begin{lemma}\label{lem:totgeod}
Let $(M,g,V,E)$ be an electrostatic system with $V^{-1}(0)=\partial M$ and $\Sigma_0\subset M$ a compact, connected, two-sided, almost properly embedded surface of finite area with $\Sigma_0\cap\partial M=\partial\Sigma_0\neq\emptyset$, which is minimal and stable for variations with compact support in $\inte\Sigma_0$. Then $\Sigma_0$ is totally geodesic, $E$ is normal to $\Sigma_0$, and $L_{\Sigma_0}V=0$.
\end{lemma}

\begin{proof}
Let $N$ be the unit normal and $A$ the second fundamental form of $\inte\Sigma_0$, and write $E=\langle E,N\rangle N+E^T$. Since $\inte\Sigma_0$ is minimal, $\Delta_{\Sigma_0}V=\Delta V-\Hess V(N,N)$, and \eqref{eq:electrostatic} gives $\Delta V=(|E|^2-\Lambda)V$ and $\Hess V(N,N)=(\Ric(N,N)-\Lambda+2\langle E,N\rangle^2-|E|^2)V$. Hence $\Delta_{\Sigma_0}V=(2|E^T|^2-\Ric(N,N))V$, that is,
\begin{equation}\label{eq:LV}
L_{\Sigma_0}V=\bigl(|A|^2+2|E^T|^2\bigr)V .
\end{equation}
Fix a non-decreasing $\chi\in C^\infty(\R)$ with $\chi=0$ on $(-\infty,1]$ and $\chi=1$ on $[2,\infty)$, and set $\eta_\varepsilon=\chi(V/\varepsilon)$ on $\Sigma_0$. Since $V$ is continuous on $\Sigma_0$, positive on $\inte\Sigma_0$ and zero on $\partial\Sigma_0$, $\eta_\varepsilon$ has compact support in $\inte\Sigma_0$. For $\varphi=\eta_\varepsilon V$, integration by parts and \eqref{eq:LV} give
\begin{align*}
0\le-\int_{\Sigma_0}\varphi\,L_{\Sigma_0}\varphi&=\int_{\Sigma_0}V^2|\nabla^{\Sigma_0}\eta_\varepsilon|^2-\int_{\Sigma_0}\eta_\varepsilon^2V\,L_{\Sigma_0}V\\
&=\int_{\Sigma_0}V^2|\nabla^{\Sigma_0}\eta_\varepsilon|^2-\int_{\Sigma_0}\eta_\varepsilon^2V^2\bigl(|A|^2+2|E^T|^2\bigr).
\end{align*}
The function $V^2|\nabla^{\Sigma_0}\eta_\varepsilon|^2=(V/\varepsilon)^2\chi'(V/\varepsilon)^2|\nabla^{\Sigma_0}V|^2$ is bounded independently of $\varepsilon$ and vanishes outside $\Sigma_0\cap\{0<V\le2\varepsilon\}$, whose area tends to zero as $\varepsilon\to0$ because $|\Sigma_0|<\infty$. Since $\eta_\varepsilon$ increases to $1$ on $\inte\Sigma_0$ as $\varepsilon\downarrow0$, monotone convergence gives $\int_{\Sigma_0}V^2(|A|^2+2|E^T|^2)\le0$. As $V>0$ on $\inte\Sigma_0$, $A\equiv0$ and $E^T\equiv0$, and then $L_{\Sigma_0}V=0$ by \eqref{eq:LV}.
\end{proof}

\begin{remark}\label{rem:almostflow}
Least area disks and annuli bounded by curves in $\partial M$ and not contained in $\partial M$ are properly embedded (Theorems~\ref{thm:plateau} and~\ref{thm:douglas} and the paragraphs after them). Such a surface $\Sigma$ is locally area-minimizing if it is disjoint from the part of the boundary of the region of minimization which is not contained in $\partial M$. Indeed, for every deformation $F$ of $\Sigma$ of small size, $F(1,\cdot)$ is a disk or an annulus, possibly immersed, in that region, with the same boundary as $\Sigma$, so $|F(1,\Sigma)|\ge|\Sigma|$. Thus Theorem~\ref{thm:A} applies to these surfaces.
\end{remark}

\subsection{Constant area along the flow}

\begin{definition}\label{def:constantarea}
Let $\Sigma\subset M$ be a compact, connected, two-sided minimal surface, either closed and contained in $\inte M$, or almost properly embedded with $\Sigma\cap\partial M=\partial\Sigma\neq\emptyset$ and bounded second fundamental form away from its exceptional set. We say that $\Sigma$ has \emph{constant area along the flow} if, for the flow $\Phi$ of Lemma~\ref{lem:flow} in the direction of one of its unit normals, $t\mapsto|\Phi(t,\Sigma)|$ is constant on some interval $[0,\delta)$.
\end{definition}

By Theorem~\ref{thm:CLS}, $\Sigma$ has constant area along the flow as soon as $|\Phi(t,\Sigma)|\ge|\Sigma|$ for small $t\ge0$. This is the case if $\Sigma$ is locally area-minimizing (Definition~\ref{def:lam}). Indeed, Lemma~\ref{lem:flow} applies to $\Sigma$, since, if $\partial\Sigma\neq\emptyset$, $\Sigma$ is totally geodesic by Lemma~\ref{lem:totgeod}; by \eqref{eq:distflow} and the estimate $|d\Phi(t,\cdot)-I|\le C_1t$ in its proof, $F(s,x)=\Phi(st,x)$ is a deformation of $\Sigma$ of size at most $Ct$, and hence $|\Phi(t,\Sigma)|\ge|\Sigma|$ for small $t$.

\subsection{Splitting along closed surfaces}

\begin{proposition}\label{prop:globalsplit}
Let $(M,g,V,E)$ be an electrostatic system with $V^{-1}(0)=\partial M$ and $\Sigma\subset\inte M$ a closed, connected, two-sided minimal surface with unit normal $N$, which has constant area along the flow in the direction $N$. Then there are a constant $E_0$ and a solution $V(s)$ of $V''=(E_0^2-\Lambda)V$, positive on a maximal open interval $I\ni0$, such that, with $\bar I$ the union of $I$ and its finite endpoints, $(\Sigma,g_\Sigma)$ has constant Gaussian curvature $\Lambda+E_0^2$, $(\bar I\times\Sigma,ds^2+g_\Sigma,V(s),E_0\partial_s)$ is a product system, and $F(s,x)=\exp_x(sN(x))$ is defined on $\bar I\times\Sigma$ and is a Riemannian covering map onto $M$ with $V\circ F=V(s)$ and $F^*E=E_0\partial_s$. Consequently $(M,g,V,E)$ is of product type, isometric to $(\bar I\times\Sigma)/\Gamma$, where $\Gamma$ is the deck group of $F$.
\end{proposition}

\begin{proof}
Let $\Phi$ be the flow of Lemma~\ref{lem:flow} in the direction $N$, with $|\Phi(t,\Sigma)|$ constant for $t\in[0,\delta)$. By Theorem~\ref{thm:CLS}(a), on $U=\Phi([0,\delta)\times\Sigma)$ the metric is $(V\circ\Phi)^2dt^2+g_\Sigma$ with $V\circ\Phi$ depending only on $t$, the leaves are totally geodesic, and $E$ is normal to them with $|E|$ constant. Reparametrizing by arclength $s=\int_0^tV(\Phi(\tau,\cdot))\,d\tau$, the curves $s\mapsto\Phi(t(s),x)$ are unit speed $g$-geodesics orthogonal to the leaves, so $\Phi(t(s),x)=F(s,x)$. Hence, for some $\varepsilon>0$ and a constant $E_0$, $F^*g=ds^2+g_\Sigma$, $V\circ F$ depends only on $s$ and $F^*E=E_0\partial_s$ on $[0,\varepsilon)\times\Sigma$, and Lemma~\ref{lem:product} shows that $V\circ F=V(s)$ solves $V''=(E_0^2-\Lambda)V$ and that $\Sigma$ has curvature $\Lambda+E_0^2$.

For $x\in\Sigma$ let $(S_-(x),S_+(x))$ be the maximal interval on which the $g$-geodesic $s\mapsto F(s,x)$ is defined in $\inte M$. Since $M$ is complete, a geodesic in $\inte M$ with $S_+(x)<\infty$ converges to a point of $\partial M$ as $s\to S_+(x)$. The domain $D=\{(s,x):S_-(x)<s<S_+(x)\}$ is open and connected, and $F:D\to\inte M$ is real-analytic by Proposition~\ref{prop:analytic}. By the first step, the analytic tensors $F^*g$ and $ds^2+g_\Sigma$, the functions $V\circ F$ and $V(s)$, and the vector fields $F^*E$ and $E_0\partial_s$ coincide on the open set $(0,\varepsilon)\times\Sigma$, hence on $D$ by Lemma~\ref{lem:continuation}. In particular $V(F(s,x))=V(s)$, so $V(s)>0$ on $(S_-(x),S_+(x))$ and $V(s)\to0$ as $s\to S_\pm(x)$ whenever $S_\pm(x)$ is finite. Since $V>0$ on $\inte M$, this shows $(S_-(x),S_+(x))=I$ for every $x$, so $D=I\times\Sigma$, and $F$ extends continuously to $\bar I\times\Sigma$, mapping $(\bar I\setminus I)\times\Sigma$ into $\partial M$. By Lemma~\ref{lem:product}, $(\bar I\times\Sigma,ds^2+g_\Sigma,V(s),E_0\partial_s)$ is a product system.

The map $F:I\times\Sigma\to\inte M$ is a local isometry. We claim it has the path lifting property. Let $c:[0,1]\to\inte M$ be a path and $\tilde c:[0,\tau)\to I\times\Sigma$ a lift. Write $\tilde c(0)=(s_0,x_0)$ and let $L$ be the length of $c$. Since $F$ is a local isometry, the lift has length at most $L$, so it stays in $\{|s-s_0|\le L\}$; since $V\ge\eta>0$ on $c([0,1])$, it also stays in $\{V(s)\ge\eta\}$. The set $\{(s,x):s\in I,\ V(s)\ge\eta,\ |s-s_0|\le L\}$ is compact, because $V(s)\to0$ at the finite endpoints of $I$. Hence $\tilde c(t)$ converges as $t\to\tau$, and the lift extends beyond $\tau$ because $F$ is a local diffeomorphism. A local diffeomorphism with the path lifting property onto a connected manifold is a covering map \cite[Chapter~7, Lemma~3.3]{dC92}; hence $F:I\times\Sigma\to\inte M$ is a Riemannian covering. Its deck group $\Gamma$ acts freely and properly discontinuously by isometries preserving $V(s)$ and $E_0\partial_s$, these isometries extend to $\bar I\times\Sigma$, and they act freely there: an isometry fixing a point $p$ of $(\bar I\setminus I)\times\Sigma$ preserves the inward normal at $p$, hence fixes the geodesic it generates, which enters $I\times\Sigma$. Since $F$ maps $(\bar I\setminus I)\times\Sigma$ into $\partial M$ and $M=\overline{\inte M}$, $F$ induces an isometry $(\bar I\times\Sigma)/\Gamma\to M$.
\end{proof}

\begin{remark}[Sub-static systems]\label{rem:substatic}
By \eqref{eq:electrostatic}, $V\Ric-\Hess V+(\Delta V)g=2V(|E|^2g-E^\flat\otimes E^\flat)\ge0$, so electrostatic systems are sub-static in the sense of \cite{CFMR26}, and \cite[Theorem~B]{CFMR26} applies to them; the warped product it gives near a closed, two-sided, embedded, locally area-minimizing surface on which $V>0$ is a product, by Theorem~\ref{thm:CLS}(a). Its global splitting assumes $\int_0^\infty V(\gamma(t))^{-1}dt=\int_0^\infty V(\gamma(t))\,dt=\infty$ for every diverging unit-speed curve $\gamma$ in the end, whereas Proposition~\ref{prop:globalsplit} uses the real-analyticity of electrostatic systems and needs no assumption on the ends. For instance, the ends $s\to+\infty$ of the systems (i)--(iii) and (i$'$)--(iii$'$) of Section~\ref{sec:models} do not satisfy the first condition, since $\int^\infty V(s)^{-1}ds<\infty$.
\end{remark}

\section{Proofs of Theorems~\ref{thm:A} and~\ref{thm:B}}\label{sec:AB}

\begin{lemma}\label{lem:surfhemisphere}
Let $(\Sigma,h)$ be a compact connected Riemannian surface with non-empty, piecewise smooth boundary and constant Gaussian curvature $K>0$, and let $V\ge0$ be a non-trivial function, of class $C^2$ up to the boundary, with
\[
\Hess_hV=-KV\,h\quad\text{in }\Sigma,\qquad V=0\quad\text{on }\partial\Sigma .
\]
Then $\partial\Sigma$ is a smooth closed geodesic, $(\Sigma,h)$ is isometric to the closed hemisphere $\Sph^2_+(K)$, and $V$ is a positive multiple of the height function.
\end{lemma}

\begin{proof}
Since $d(|\nabla V|^2+KV^2)=2\Hess_hV(\nabla V,\cdot)+2KV\,dV=0$ and $V\not\equiv0$, we have $|\nabla V|^2+KV^2=c^2$ for a constant $c>0$, so $|\nabla V|=c$ wherever $V=0$. An interior zero, a minimum of $V\ge0$, or a corner of $\partial\Sigma$, where $V$ vanishes along two independent directions, would have $\nabla V=0$; hence $V>0$ in $\inte\Sigma$ and $\partial\Sigma$ is a regular level set of $V$, with geodesic curvature $\pm\Hess_hV(T,T)/c=0$ for a unit tangent $T$. By the Gauss--Bonnet theorem, $K|\Sigma|=2\pi\chi(\Sigma)$, so $\chi(\Sigma)>0$; as $\partial\Sigma\neq\emptyset$, $\Sigma$ is a disk and $|\Sigma|=2\pi/K$. In Fermi coordinates along the geodesic $\partial\Sigma$, the metric is $ds^2+\cos^2(\sqrt K s)\,d\theta^2$, which is even in $s$; hence the double of $\Sigma$ along $\partial\Sigma$ is a smooth closed surface of constant curvature $K$, diffeomorphic to $\Sph^2$, and therefore isometric to the round sphere $\Sph^2(K)$. In it, $\partial\Sigma$ is a closed geodesic, i.e.\ a great circle, and $\Sigma$ is one of the two hemispheres it bounds. Finally, the solutions of $\Hess f=-Kf\,h$ on $\Sph^2(K)\subset\R^3$ are the restrictions of linear functions, since such a solution is determined by its value and gradient at one point; vanishing on $\partial\Sigma$ and positive in $\inte\Sigma$, $V$ is a positive multiple of the height function.
\end{proof}

\begin{lemma}\label{lem:hemisphere}
Let $(M,g,V,E)$ be an electrostatic system with $V^{-1}(0)=\partial M\neq\emptyset$ such that $E=0$ and $g$ has constant sectional curvature $K>0$ on a non-empty open subset of $\inte M$. Then $K=\Lambda/3$, $E\equiv0$, $(M,g)$ is isometric to the closed hemisphere $\Sph^3_+(K)$, and $V$ is a positive multiple of the height function; that is, $(M,g,V,E)$ is the de Sitter system.
\end{lemma}

For compact $M$, once $E\equiv0$ and the curvature is constant on all of $M$, the conclusion is the version of Obata's theorem \cite{Oba62} for manifolds with boundary, due to Reilly \cite{Rei77}; for static vacuum metrics it also follows from the classification of locally conformally flat static triples \cite{Kob82,Laf83} (see \cite[Theorem~1]{Amb17}). In the proof of Theorem~\ref{thm:A}, these conditions are known only on an open set and $M$ is not known to be compact; the proof below reduces to that case.

\begin{proof}
The function $|\Riem-\tfrac K2\,g\owedge g|^2+|E|^2$ is real-analytic on the connected manifold $\inte M$ by Proposition~\ref{prop:analytic} and vanishes on an open set, so it vanishes on $M$ by Lemma~\ref{lem:continuation} and continuity. Thus $E\equiv0$ and $\Ric=2Kg$; Lemma~\ref{lem:basic}(i) gives $6K=R_g=2\Lambda$, so $K=\Lambda/3$, and the first equation in \eqref{eq:electrostatic} gives $\Hess V=V({\Ric}-\Lambda g)=-KVg$. Since $\nabla(|\nabla V|^2+KV^2)=2\Hess V(\nabla V)+2KV\nabla V=0$, there is $c>0$ with $|\nabla V|^2+KV^2=c^2$ on $M$. Let $x\in\partial M$ and let $\sigma$ be the inward normal unit speed geodesic from $x$. Then $y=V\circ\sigma$ solves $y''=-Ky$, $y(0)=0$, $y'(0)=c$, so $y(s)=(c/\sqrt K)\sin(\sqrt Ks)>0$ for $0<s<\pi/\sqrt K$; hence $\sigma$ does not reach $\partial M$ before $s=\pi/\sqrt K$, and $p=\sigma(\pi/(2\sqrt K))$ satisfies $V(p)=c/\sqrt K$ and $\nabla V(p)=0$. For every unit $v\in T_pM$, $V(\exp_p(sv))=(c/\sqrt K)\cos(\sqrt Ks)$, so $\exp_p$ is defined on the closed ball $\bar B$ of radius $\pi/(2\sqrt K)$, maps its interior into $\inte M$ and its boundary sphere into $\partial M$, and is a local isometry from $\bar B$ with the metric of $\Sph^3_+(K)$. Its image is compact, and open because $\exp_p$ is a local isometry of manifolds with boundary; hence $M=\exp_p(\bar B)$ is compact. By Reilly's theorem \cite{Rei77}, $(M,g)$ is isometric to $\Sph^3_+(K)$, and $V$, which solves $\Hess V=-KVg$ and vanishes on $\partial M$, is a positive multiple of the height function, as in the proof of Lemma~\ref{lem:surfhemisphere}.
\end{proof}

\begin{proof}[Proof of Theorem~\ref{thm:A}]
By Definition~\ref{def:constantarea} and the paragraph after it, the flow $\Phi$ of Lemma~\ref{lem:flow} starting at $\Sigma$, defined on $\Sigma\setminus S$ with $S\subset\partial\Sigma$ the finite exceptional set, has constant area on some interval $[0,\delta)$; only this is used (Remark~\ref{rem:constantarea}). By Theorem~\ref{thm:CLS}(b), $\Lambda>0$; $E=0$ on $U=\Phi((0,\delta)\times\inte\Sigma)$, which is open in $\inte M$ because $\partial_t\Phi=(V\circ\Phi)N_\Phi$ is transverse to the leaves; the leaves are totally geodesic and $\Phi^*g=(V\circ\Phi)^2dt^2+g_\Sigma$; and $(\Sigma,g_\Sigma)$ is a domain bounded by geodesic arcs in the round sphere of curvature $K=\Lambda/3$, on which $V\circ\Phi(t,\cdot)$ is, for each $t$, a static potential vanishing on $\partial\Sigma$; as in the proof of Lemma~\ref{lem:surfhemisphere}, such a potential is the restriction of a linear function on $\R^3\supset\Sph^2(K)$, hence smooth up to $\partial\Sigma$.

\emph{The leaves are hemispheres.} By Lemma~\ref{lem:surfhemisphere}, $(\Sigma,g_\Sigma)$ is isometric to the hemisphere $\Sph^2_+(K)$, and, for each $t$, $V\circ\Phi(t,\cdot)$ is a positive multiple of its height function. In spherical coordinates $(\theta,\phi)\in[0,\frac\pi2]\times[0,2\pi)$ on $\Sigma$ centered at the pole, this means that
\begin{equation}\label{eq:pullbackA}
V(\Phi(t,\theta,\phi))=a(t)\cos\theta,\qquad \Phi^*g=a(t)^2\cos^2\theta\,dt^2+\frac3\Lambda\bigl(d\theta^2+\sin^2\theta\,d\phi^2\bigr),
\end{equation}
for a smooth positive function $a$.

\emph{The function $a$.} Since $E=0$ on $U$, the second equation in~\eqref{eq:electrostatic} reads $\Delta V=-\Lambda V$ there. For the metric~\eqref{eq:pullbackA},
\[
-\Lambda a\cos\theta=\Delta V=\frac{a\ddot a-\dot a^2}{a^3\cos\theta}-\frac\Lambda3\,a\,\frac{2\cos^2\theta-\sin^2\theta}{\cos\theta}\implies a\ddot a-\dot a^2+\frac\Lambda3\,a^4=0 .
\]
The positive solutions of this equation are
\[
a(t)=\frac{a_0}{\cosh\bigl(a_0\sqrt{\Lambda/3}\,(t-t_0)\bigr)},\qquad a_0>0,\ t_0\in\R .
\]

\emph{Constant curvature.} Let $r(t)=\arctan\sinh\bigl(a_0\sqrt{\Lambda/3}\,(t-t_0)\bigr)$. Then $\dot r=\sqrt{\Lambda/3}\,a$, so $a^2dt^2=\frac3\Lambda\,dr^2$, and~\eqref{eq:pullbackA} becomes
\[
\Phi^*g=\frac3\Lambda\bigl(\cos^2\theta\,dr^2+d\theta^2+\sin^2\theta\,d\phi^2\bigr),
\]
which is the round metric of curvature $\Lambda/3$ (Remark~\ref{rem:openbook}). Thus $E=0$ and $g$ has constant sectional curvature $\Lambda/3$ on the open set $U$, and, by Lemma~\ref{lem:hemisphere}, $(M,g,V,E)$ is the de Sitter system $\Sph^3_+(K)$.

Finally, $\Sigma$ is totally geodesic and $\Sigma\setminus S$ is connected, so $\Sigma$ lies in a great sphere $P$ of $\Sph^3(K)$, different from $\partial\Sph^3_+(K)$ because $\Sigma$ meets $\inte M$. Then $\inte\Sigma$ is open and closed in the interior of the hemisphere $H=P\cap\Sph^3_+(K)$, which is bounded by the great circle $P\cap\partial\Sph^3_+(K)$, so $\Sigma=H$. The flow starting at $H$ is described in Remark~\ref{rem:openbook}.
\end{proof}

\begin{remark}\label{rem:openbook}
Let $K=\Lambda/3$ and, for $r\in(-\frac\pi2,\frac\pi2)$, $\theta\in[0,\frac\pi2]$ and $\phi\in[0,2\pi)$, let
\[
R=\arccos(\cos\theta\cos r)\in\bigl[0,\tfrac\pi2\bigr],\qquad
\vartheta=\arccos\frac{\cos\theta\sin r}{\sqrt{1-\cos^2\theta\cos^2r}}\in[0,\pi],\qquad
\varphi=\phi .
\]
Then
\[
\frac1K\bigl(\cos^2\theta\,dr^2+d\theta^2+\sin^2\theta\,d\phi^2\bigr)=\frac1K\Bigl(dR^2+\sin^2R\,\bigl(d\vartheta^2+\sin^2\vartheta\,d\varphi^2\bigr)\Bigr),
\]
the round metric of the hemisphere $\Sph^3_+(K)$ in polar coordinates about its pole. Indeed, both are induced by
\[
x=K^{-1/2}\bigl(\sin\theta\cos\phi,\ \sin\theta\sin\phi,\ \cos\theta\sin r,\ \cos\theta\cos r\bigr)\in\Sph^3(K)\subset\R^4,
\]
for which
\[
K^{1/2}x=\bigl(\sin R\sin\vartheta\cos\varphi,\ \sin R\sin\vartheta\sin\varphi,\ \sin R\cos\vartheta,\ \cos R\bigr),
\]
so $x_4\ge0$. In the proof of Theorem~\ref{thm:A}, $a=a_0\cos r$, so $V=a_0\cos r\cos\theta=a_0\cos R$ is a multiple of the height function $x_4$. The boundary of each leaf $\Sigma_t=\{r=r(t)\}$ is $\{\theta=\frac\pi2\}$, which corresponds to $(R,\vartheta)=(\frac\pi2,\frac\pi2)$, that is, to the great circle $\{x_3=x_4=0\}$ of $\partial\Sph^3_+(K)=\{x_4=0\}$. Since each $\Sigma_t$ is totally geodesic, it is the intersection of $\Sph^3_+(K)$ with the hyperplane $\{x_3=(\tan r(t))\,x_4\}$, which contains this circle; as $t\to\pm\infty$, $r(t)\to\pm\frac\pi2$, and $\Sigma_t$ converges to the two hemispheres $\{x_4=0,\ \pm x_3\ge0\}$ of $\partial\Sph^3_+(K)$. Therefore $\{\Sigma_t\}_{t\in\R}$ is the open book foliation of $\Sph^3_+(K)$ with binding $\{x_3=x_4=0\}$.
\end{remark}

\begin{proof}[Proof of Theorem~\ref{thm:B}]
By Remark~\ref{rem:constantarea}, we only use that $\Sigma$, or, if $\Sigma$ is one-sided, its two-sided preimage in the double cover of Definition~\ref{def:lam}(ii), has constant area along the flow.

\emph{Step 1: product structure.} If $\Sigma$ is two-sided, Proposition~\ref{prop:globalsplit} applies with $\hat\Sigma=\Sigma$. If $\Sigma$ is one-sided, let $\pi:\hat M\to M$ be the double cover of Definition~\ref{def:lam}(ii); Proposition~\ref{prop:globalsplit} applied to $\hat\Sigma=\pi^{-1}(\Sigma)$ in $\hat M$ gives $\hat M\cong(\bar I\times\hat\Sigma)/\hat\Gamma$, and $M$ is the quotient of $\bar I\times\hat\Sigma$ by the group generated by $\hat\Gamma$ and a lift of the deck involution of $\hat M$. In both cases we obtain a product system
\[
X=\bar I\times\hat\Sigma,\qquad g_X=ds^2+h,\qquad V=V(s),\qquad E=E_0\partial_s,
\]
in which $(\hat\Sigma,h)$ has constant curvature $K=\Lambda+E_0^2$ and $V''=\mu V$ with $\mu=E_0^2-\Lambda$, together with a Riemannian covering $F:X\to M$ such that $V\circ F=V(s)$ and $F^*E=E_0\partial_s$, and $M\cong X/\Gamma$, where $\Gamma$ is the deck group of $F$. We record the relations
\begin{equation}\label{eq:BKmu}
K+\mu=2E_0^2,\qquad K-\mu=2\Lambda .
\end{equation}
Since $F\circ\gamma=F$, every $\gamma\in\Gamma$ is an isometry of $g_X$ which preserves the orientation of $X$ induced from $M$ and satisfies
\begin{equation}\label{eq:Bdeck}
V\circ\gamma=V,\qquad d\gamma(E_0\partial_s)=E_0\partial_s .
\end{equation}
Finally, $|E|=|E_0|$, so Lemma~\ref{lem:basic}(i) gives $R_g=2K$ on $M$. Thus $K>0$ exactly when $R_g>0$ somewhere on $\Sigma$, and in particular when $\Lambda>0$.

\emph{Step 2: deck transformations are products.} We claim that every $\gamma\in\Gamma$ has the form
\begin{equation}\label{eq:Bform}
\gamma(s,x)=\bigl(\sigma_\gamma(s),\varphi_\gamma(x)\bigr),
\end{equation}
where $\sigma_\gamma$ is an isometry of $\bar I$ with $V\circ\sigma_\gamma=V$ and $\varphi_\gamma$ is an isometry of $(\hat\Sigma,h)$. Suppose first that $V$ is not constant. Since $\gamma$ is an isometry preserving $V$, we have $d\gamma(\nabla V)=\nabla V\circ\gamma$, and $\nabla V=V'(s)\,\partial_s$. The function $V'$ is a non-trivial solution of $y''=\mu y$, so its zeros are isolated; hence $d\gamma(\partial_s)=\pm\partial_s$ on a dense open set and, by continuity and connectedness,
\begin{equation}\label{eq:Beps}
d\gamma(\partial_s)=\epsilon\,\partial_s\qquad\text{for a constant }\epsilon\in\{1,-1\}.
\end{equation}
If $E_0\neq0$, then \eqref{eq:Beps} holds with $\epsilon=1$, by \eqref{eq:Bdeck}. Now \eqref{eq:Beps} implies \eqref{eq:Bform}: the function $s\circ\gamma$ has gradient $\epsilon\,\partial_s$, so $s\circ\gamma=\epsilon s+c$ for a constant $c$, and the projection $p:X\to\hat\Sigma$ satisfies $d(p\circ\gamma)(\partial_s)=0$, so $p\circ\gamma$ does not depend on $s$. In the remaining case, $V$ is constant and $E_0=0$. Then $\mu=0$ and, by \eqref{eq:BKmu}, $K=\Lambda=0$, so $\hat\Sigma$ is a flat torus and $X=\R\times\hat\Sigma$ is flat. An isometry of $X$ lifts to a Euclidean isometry of $\R^3$ normalizing the lattice of translations of $\hat\Sigma$, hence preserving its span and the orthogonal line, which gives \eqref{eq:Bform} again.

\emph{Step 3: orientation and electric field.} An isometry of $\bar I$ is a translation or a reflection. By \eqref{eq:Bform}, $\gamma$ preserves the orientation of $X$ if and only if $\sigma_\gamma$ and $\varphi_\gamma$ both preserve or both reverse the orientation; hence
\begin{equation}\label{eq:Borient}
\varphi_\gamma\text{ preserves the orientation of }\hat\Sigma\iff\sigma_\gamma\text{ is a translation}.
\end{equation}
A reflection maps $E_0\partial_s$ to $-E_0\partial_s$, so, by \eqref{eq:Bdeck},
\begin{equation}\label{eq:Brefl}
\sigma_\gamma\text{ is a reflection}\implies E_0=0 .
\end{equation}

\emph{Step 4: reduction to a group acting on $\bar I$.} By \eqref{eq:Bform}, $\gamma\mapsto\sigma_\gamma$ is a homomorphism from $\Gamma$ to the isometry group of $\bar I$. Its kernel $\Gamma_0$ consists of the elements of the form $(s,x)\mapsto(s,\varphi(x))$. If such a $\varphi$ fixed a point $x$, the element would fix $(s,x)$; as $\Gamma$ acts freely, $\Gamma_0$ acts freely on $\hat\Sigma$, by orientation-preserving isometries by \eqref{eq:Borient}. Hence $\hat\Sigma/\Gamma_0$ is a closed orientable surface of constant curvature $K$, the cross-section of the model, which we denote again by $\Sigma$, as in Section~\ref{sec:models}. This agrees with the notation of the statement when the given surface $\Sigma$ is two-sided: then $\hat\Sigma=\Sigma$ and $F(0,x)=x$, so an element $(s,x)\mapsto(s,\varphi(x))$ of $\Gamma_0$ satisfies $\varphi(x)=F(0,\varphi(x))=F(0,x)=x$, and $\Gamma_0$ is trivial. If $\Sigma$ is one-sided, the cross-section is a different surface; for the quotient of the Nariai system, $\Sigma$ is a projective plane and the cross-section is a sphere. Thus $X/\Gamma_0=\bar I\times\Sigma$, and, since $\Gamma_0$ is normal,
\[
M\cong(\bar I\times\Sigma)/G,\qquad G=\Gamma/\Gamma_0,
\]
where $G$ acts faithfully on $\bar I$ through the isometries $\sigma_\gamma$. It remains to determine $\Sigma$ and $G$.

\emph{Step 5: isometries of $\bar I$ preserving $V$.} After a translation and possibly a reflection of $s$ (changing the sign of $E_0$), $V$ is one of the maximal positive solutions of $V''=\mu V$ listed in Section~\ref{sec:models}. The isometries of $\R$ are the translations and the reflections, the only isometry of $[0,\infty)$ is the identity, and those of $[-a,a]$ are the identity and $s\mapsto-s$. Checking each solution, we find that a non-trivial translation preserves $V$ only if $V\equiv1$, which forces $\mu=0$, and that a reflection preserves $V$ only if $V\equiv1$, or $V$ is $\cosh(\sqrt\mu\,s)$ or $\cos(\sqrt{-\mu}\,s)$ and the reflection is $s\mapsto-s$. With \eqref{eq:Brefl}, we conclude that $G$ is trivial unless $V\equiv1$, or $E_0=0$ and $V$ is $\cosh(\sqrt\mu\,s)$ or $\cos(\sqrt{-\mu}\,s)$.

\emph{Step 6: spherical cross-sections.} Let $K>0$. Then $\hat\Sigma$ is the round sphere of curvature $K$, with isometry group $O(3)$. An element of $\Gamma_0$ acts on it by a rotation, which has a fixed point; since $\Gamma_0$ acts freely, it is trivial. Thus $\Sigma=\Sph^2$, $G=\Gamma$, and $\gamma\mapsto\sigma_\gamma$ is injective. By Step~5, there are three cases.

(a) \emph{$\Gamma$ is trivial.} Then $M=\bar I\times\Sph^2$, and the six solutions give the systems (i)--(iv), (vi) and (vii) of Table~\ref{tab:spherical}.

(b) \emph{$V\equiv1$ and $\Gamma$ is non-trivial.} Then $\mu=0$, so $E_0^2=\Lambda=K/2>0$ by \eqref{eq:BKmu}, and $\Gamma$ consists of translations by \eqref{eq:Brefl}. It is discrete, because it acts properly discontinuously on $\R\times\Sph^2$ and $\Sph^2$ is compact. Hence $\Gamma$ is generated by some $(s,p)\mapsto(s+L,Ap)$ with $L>0$, where $A\in SO(3)$ by \eqref{eq:Borient}. This is (v).

(c) \emph{$\Gamma$ contains an element $\gamma$ with $\sigma_\gamma$ a reflection.} Then $E_0=0$ by \eqref{eq:Brefl}, so $K=\Lambda>0$ and $\mu=-\Lambda<0$ by \eqref{eq:BKmu}; hence $V=\cos(\sqrt\Lambda\,s)$ on $[-a,a]$ and $\sigma_\gamma(s)=-s$. By injectivity, $\Gamma=\{\mathrm{id},\gamma\}$ and $\gamma^2=\mathrm{id}$, so $\varphi=\varphi_\gamma$ satisfies $\varphi^2=\mathrm{Id}$, and $\varphi\in O(3)\setminus SO(3)$ by \eqref{eq:Borient}. Write $\varphi=-R$ with $R\in SO(3)$. Then $R^2=\mathrm{Id}$, so $R$ is the identity or a rotation by $\pi$. Since $\gamma$ maps $\{0\}\times\Sph^2$ to itself by $\varphi$ and acts freely, $-R$ has no fixed point on $\Sph^2$, that is, $-1$ is not an eigenvalue of $R$; this excludes the rotation by $\pi$. Hence $\varphi=-\mathrm{Id}$ and $\gamma(s,p)=(-s,-p)$, which is (viii).

\emph{Step 7: flat and hyperbolic cross-sections.} Let $K\le0$. By the Gauss--Bonnet theorem, $\Sigma$ has positive genus: it is a flat torus if $K=0$ and a hyperbolic surface if $K<0$. By \eqref{eq:BKmu}, $\Lambda=K-E_0^2\le0$ and $\mu=2E_0^2-K\ge0$, so $V$ is not $\cos(\sqrt{-\mu}\,s)$.

Suppose that $\mu>0$. By Step~5, $G$ is trivial unless $E_0=0$ and $V=\cosh(\sqrt\mu\,s)$. If $G$ is trivial, then $M=\bar I\times\Sigma$, and we obtain (i$'$)--(iii$'$) of Table~\ref{tab:flathyp}. Otherwise $E_0=0$, so $\Lambda=K$ and $\mu=-K>0$; thus $\Sigma$ is hyperbolic and $G=\{\mathrm{id},\gamma\}$, with $\gamma(s,x)=(-s,\iota(x))$ for an isometry $\iota$ of $\Sigma$. Here $\iota^2=\mathrm{id}$, because $\gamma^2$ lies in $\Gamma_0$ and so acts trivially on $\Sigma=\hat\Sigma/\Gamma_0$; $\iota$ reverses the orientation, by \eqref{eq:Borient}; and $\iota$ has no fixed points, because $\gamma$ acts freely on $\{0\}\times\Sigma$. This is (viii$'$).

Suppose that $\mu=0$. By \eqref{eq:BKmu}, $E_0=0$ and $K=\Lambda=0$, so $\Sigma$ is a flat torus. If $V=s$ on $[0,\infty)$, then $G$ is trivial and we obtain (vi$'$). If $V\equiv1$, every isometry of $\R$ preserves $V$, and $G$ is a discrete group of isometries of $\R$ with no further restriction; these are the flat systems (iv$'$).
\end{proof}

\section{Compact electrostatic systems}\label{sec:compact}

We first prove Theorem~\ref{thm:D}, on which Theorem~\ref{thm:C} and the corollaries depend.

\begin{proof}[Proof of Theorem~\ref{thm:D}]
(i) Since $\Lambda+|E|^2>0$ by Lemma~\ref{lem:basic}(v), every component of $\partial M$ is a sphere by \cite[Theorem~A]{CLS24}. In particular $H_1(\partial M;\R)=0$. Suppose $H_1(M;\R)\neq0$. By Lefschetz duality $H_2(M,\partial M;\R)\cong H^1(M;\R)\neq0$, and the exact sequence
\[
H_2(\partial M;\R)\xrightarrow{\ i_*\ }H_2(M;\R)\longrightarrow H_2(M,\partial M;\R)\longrightarrow H_1(\partial M;\R)=0
\]
yields $\alpha\in H_2(M;\Z)$ whose image in $H_2(M;\R)$ does not lie in $i_*H_2(\partial M;\R)$. Since $\partial M$ is totally geodesic, $M$ is mean-convex, and Theorem~\ref{thm:homology} gives a minimizer $\sum n_i[\Sigma_i]$ of $\alpha$. Its components lying on $\partial M$ are components of $\partial M$, so some $\Sigma_i$ is a closed, embedded, orientable, locally area-minimizing surface in $\inte M$. As $\Lambda>0$, Theorem~\ref{thm:B} shows that $(M,g,V,E)$ is one of the systems (i)--(viii); being compact with non-empty boundary, it is (vii) or (viii), both of which have $H_1(M;\R)=0$, a contradiction. Hence $H_1(M;\R)=0$, the map $H_2(M;\R)\to H_2(M,\partial M;\R)$ vanishes, $H_2(M;\R)=i_*H_2(\partial M;\R)$, and Lemma~\ref{lem:boundaryclasses} gives $H_2(M;\R)\cong\R^{\ell}/\R(1,\dots,1)\cong\R^{\ell-1}$.

(ii) This follows from (i) and Lemma~\ref{lem:basic}(iv).

(iii) If $\ell=1$, no non-empty union of horizons satisfies the hypothesis, so we assume $\ell\ge2$. Since $\Lambda>0$, Theorem~\ref{thm:B} shows that $\inte M$ contains no closed, embedded, orientable, locally area-minimizing surface unless $(M,g,V,E)$ is one of the systems (vii) and (viii), the only compact models with non-empty boundary; as (viii) has a single horizon, this leaves the charged Nariai systems. Assume first that $(M,g,V,E)$ is not a charged Nariai system; the argument follows the homological minimization in the proofs of \cite[Theorem~B]{Amb17} and \cite[Theorem~10]{CLS24}. By Theorem~\ref{thm:homology}, the class $\alpha=[\partial_JM]\in H_2(M;\Z)$ has a minimizer of the form $\sum_ic_i[\partial_iM]$, $c_i\in\Z$, with mass $\sum_i|c_i||\partial_iM|$. By Lemma~\ref{lem:boundaryclasses}, $c_i=\mathbf 1_J(i)+t$ for some $t\in\Z$, so the mass is
\[
f(t)=|1+t|\,|\partial_JM|+|t|\,|\partial M\setminus\partial_JM| .
\]
This function is convex and piecewise linear, with $f(0)=|\partial_JM|\le f(-1)=|\partial M\setminus\partial_JM|$, so the minimal mass in $\alpha$ is $|\partial_JM|$: every closed embedded orientable surface $\tilde\Sigma$ homologous to $\partial_JM$ satisfies $|\tilde\Sigma|\ge|\partial_JM|$, and if equality holds, $\tilde\Sigma$ is a minimizer. Its components in $\inte M$ would be closed locally area-minimizing surfaces, which we have excluded, and its components meeting $\partial M$ are horizons by Lemma~\ref{lem:maxprinciple}; hence $\tilde\Sigma=\partial_JM$, or $\tilde\Sigma=\partial M\setminus\partial_JM$ when $f(0)=f(-1)$. For a charged Nariai system, the two horizons have the same area, and a closed embedded orientable surface homologous to a horizon has area at least that of the horizon, because the projection onto the $\Sph^2$ factor does not increase area and has degree $\pm1$ on it. Finally, a horizon $\partial_jM$ with $|\partial_jM|\le|\partial M\setminus\partial_jM|$ minimizes area in its homology class, and is therefore stable. The last assertion of (iii) follows by contraposition.

(iv) The uniqueness of the unstable horizon also follows from (iii), since two unstable horizons would each have more than half of the total area. The topology of $M$ follows from \cite[Theorems~A and~10]{CLS24}, which apply because $\Lambda>0$: by \eqref{eq:integral}, $\int_M(\Lambda-|E|^2)V\,dv_g>0$ with $V>0$ in $\inte M$, so $\Lambda>\min_M|E|^2\ge0$. In particular $M$ is simply connected. Let $S=\partial_\ell M$ be the unstable horizon and suppose that its Morse index is at least two. Let $\lambda_1<\lambda_2<0$ be the first two eigenvalues of $L_S$, where $L_S=\Delta_S+\Ric(\nu,\nu)$ since $S$ is totally geodesic, and let $u_2$ be an eigenfunction for $\lambda_2$. Since $u_2$ is orthogonal to the positive first eigenfunction, it changes sign, and by Courant's nodal domain theorem \cite{CH53} it has exactly two nodal domains. By Lemma~\ref{lem:nodal}, $\Gamma=u_2^{-1}(0)$ is a smooth Jordan curve splitting $S$ into two disks $\Omega_1,\Omega_2$, and $v_i=|u_2|\,|_{\Omega_i}$ is a first Dirichlet eigenfunction of $L_S$ on $\Omega_i$ with eigenvalue $\lambda_2<0$.

Since $M$ is simply connected, $\Gamma$ is null-homotopic in $M$, and Theorem~\ref{thm:plateau} applied to $W=M$ gives a least area disk $D$ bounded by $\Gamma$ which is almost properly embedded or contained in $\partial M$. In the second case $D\subset S$, so $D=\Omega_1$ or $D=\Omega_2$. But the disks $\Omega_i^\varepsilon=\{\exp_x(\varepsilon v_i(x)\nu(x)):x\in\Omega_i\}\subset M$ are bounded by $\Gamma$ and, since $\Omega_i$ is totally geodesic,
\[
|\Omega_i^\varepsilon|=|\Omega_i|+\frac{\varepsilon^2}{2}\,\lambda_2\int_{\Omega_i}v_i^2+o(\varepsilon^2)<|\Omega_i|
\]
for small $\varepsilon>0$, which contradicts the minimality of $D$. Hence $D$ is an almost properly embedded disk bounded by $\Gamma$, of least area among such disks in $M$, so it is locally area-minimizing (Remark~\ref{rem:almostflow}), and Theorem~\ref{thm:A} shows that $(M,g,V,E)$ is the de Sitter system. Its boundary is a round sphere of curvature $K=\Lambda/3$ with $\Ric(\nu,\nu)=2K$, so $L_S=\Delta_S+2K$ has exactly one negative eigenvalue, contradicting the assumption. Hence the Morse index of $S$ is one.

(v) We use the topological classification of Carlotto and Li \cite[Theorem~1.1]{CL24} of compact $3$-manifolds carrying a metric of positive scalar curvature with strictly mean-convex boundary. Here $R_g>0$, but $\partial M$ is totally geodesic, so we first deform $g$ conformally (compare \cite[Section~3]{CL24b}). Let $\rho$ be a smooth function on $M$ which coincides with the distance to $\partial M$ near $\partial M$. For small $\varepsilon>0$, the metric $\tilde g=e^{-2\varepsilon\rho}g$ still has positive scalar curvature, since $M$ is compact, and $\partial M$ has mean curvature $2\varepsilon>0$ with respect to $\tilde g$ and the outward normal. By \cite[Theorem~1.1]{CL24}, $M$ is diffeomorphic to a connected sum of handlebodies, spherical space forms and copies of $\Sph^2\times\Sph^1$. Since the components of $\partial M$ are spheres, the handlebodies are balls, and since $b_1(M)=0$ by (i), there are no summands $\Sph^2\times\Sph^1$. Discarding summands $\Sph^3$, we conclude that $M$ is the connected sum of $k\ge0$ spherical space forms $\Sph^3/G_1,\dots,\Sph^3/G_k$ with finite non-trivial groups $G_i$ (the $3$-sphere if $k=0$), with $\ell$ open balls removed, so $\pi_1(M)=G_1*\dots*G_k$. Suppose that $k\ge2$. The kernel $N$ of the natural homomorphism $G_1*\dots*G_k\to G_1\times\dots\times G_k$ has finite index and meets every conjugate of every $G_i$ trivially, so it is free by the Kurosh subgroup theorem \cite{Kur34}; it is non-trivial, because for non-trivial $a\in G_1$ and $b\in G_2$ the element $ab$ has infinite order and some power of it lies in $N$. The finite covering of $M$ associated with $N$, with the pulled back data, is a compact electrostatic system whose first Betti number is the rank of $N$, which is positive; this contradicts (i) applied to that covering. Hence $k\le1$, which proves (v); if some horizon is unstable, $G$ is trivial by (iv).
\end{proof}

\begin{remark}\label{rem:CLSindex}
In \cite[Theorems~B and~C]{CLS24}, index one is obtained by min-max methods for an unstable surface $\Sigma=\partial\Omega_1\cap\partial\Omega_2\subset V^{-1}(0)$ separating two compact maximal regions $\Omega_1,\Omega_2$ on which $V\neq0$, under stability assumptions on the remaining components of $\partial\Omega_i$. Theorem~\ref{thm:D}(iv) concerns instead the boundary of a single compact region, with no assumption on the other horizons, and its proof uses Theorem~\ref{thm:A} rather than min-max theory.
\end{remark}

\begin{lemma}\label{lem:potential}
Let $(M,g,V,E)$ be as in Theorem~\ref{thm:D}, let $\psi_j$ be the value of $\psi$ on $\partial_jM$, and let $Q_j=\frac1{4\pi}\int_{\partial_jM}\langle E,\nu\rangle$. Then $\sum_jQ_j=0$ and
\[
\int_MV|E|^2\,dv_g=-4\pi\sum_{j=1}^{\ell}\psi_jQ_j .
\]
In particular, $E\equiv0$ if and only if $\psi$ takes the same value on all horizons, if and only if every horizon has charge zero.
\end{lemma}

\begin{proof}
Since $\dv E=0$, the divergence theorem gives $\sum_jQ_j=0$, and by Lemma~\ref{lem:basic}(iv), $\int_MV|E|^2=\int_M\dv(\psi E)=-\sum_j\psi_j\int_{\partial_jM}\langle E,\nu\rangle$. If all $\psi_j$ are equal, the right-hand side vanishes because $\sum_jQ_j=0$, and it also vanishes if all $Q_j$ vanish; conversely, $E=0$ gives $d\psi=0$ and $Q_j=0$.
\end{proof}

\begin{remark}[Smarr-type identity]\label{rem:smarr}
Combining \eqref{eq:integral} with Lemma~\ref{lem:potential} gives
\[
\Lambda\int_MV\,dv_g=\sum_{j=1}^{\ell}\kappa_j|\partial_jM|-4\pi\sum_{j=1}^{\ell}\psi_jQ_j ,
\]
which does not depend on the normalization of $\psi$, since $\sum_jQ_j=0$. It is an integrated form of the Komar identity \cite{Kom59}, and an analogue, for static systems with cosmological constant, of Smarr's mass formula \cite{Sma73}. For two horizons, the last term equals $-4\pi(\psi_1-\psi_2)Q_1=4\pi|\psi_1-\psi_2|\,|Q_1|$, since $Q_2=-Q_1$ and, by Lemma~\ref{lem:potential}, $(\psi_1-\psi_2)Q_1\le0$.
\end{remark}

\begin{proof}[Proof of Theorem~\ref{thm:masscharge}]
Since $\ell\ge2$, a horizon $S'$ of largest area satisfies $|\partial_jM|\le|S'|\le|\partial M|-|\partial_jM|$ for every other horizon $\partial_jM$, so the hypothesis on $S$ holds for every horizon other than $S'$. By Theorem~\ref{thm:D}(iii), $S$ minimizes area in its homology class, so it is stable; since $S$ is a sphere, Proposition~\ref{prop:rigidity} gives $\Lambda|S|+16\pi^2Q_S^2/|S|\le4\pi$ (compare \cite{BBB23}), that is,
\[
\Lambda r_S^4-r_S^2+Q_S^2\le0 .
\]
Hence $1-4\Lambda Q_S^2\ge0$ and $r_-\le r_S\le r_+$, since the roots of $\Lambda r^4-r^2+Q_S^2$, as a polynomial in $r^2$, are $r_\pm^2$. If $r_S\in\{r_-,r_+\}$, equality holds, Proposition~\ref{prop:rigidity} shows that the system is of product type, and, by the proof of Theorem~\ref{thm:B}, a compact system of product type with at least two horizons is a charged Nariai system. For the latter, $r_S^2=1/K$ and $Q_S^2=E_0^2/K^2$, with $K=\Lambda+E_0^2$ and $E_0^2<\Lambda$, so $1-4\Lambda Q_S^2=(\Lambda-E_0^2)^2/K^2>0$ and $r_+^2=1/K=r_S^2$. Therefore $4\Lambda Q_S^2<1$, $r_-<r_S\le r_+$, and $r_S=r_+$ if and only if the system is a charged Nariai system.

The function $\mathfrak m(r)=\frac r2\bigl(1+Q_S^2/r^2-\Lambda r^2/3\bigr)$, which is defined on $[r_-,r_+]$ (with $\mathfrak m(0)=0$ if $Q_S=0$), satisfies
\[
\mathfrak m'(r)=\frac12\Bigl(1-\frac{Q_S^2}{r^2}-\Lambda r^2\Bigr)=-\frac{\Lambda r^4-r^2+Q_S^2}{2r^2}>0\qquad\text{on }(r_-,r_+).
\]
Hence $\mathfrak m$ is strictly increasing on $[r_-,r_+]$, and, since $Q_S^2=r_\pm^2(1-\Lambda r_\pm^2)$, $\mathfrak m(r_\pm)=r_\pm(1-\frac23\Lambda r_\pm^2)=m_\pm$. As $m_H(S)=\mathfrak m(r_S)$, this gives the bounds on the charged Hawking mass and the equality case.

Finally, Lemma~\ref{lem:potential} gives $\sum_jQ_j=0$ and, if $E\not\equiv0$, $\sum_j\psi_jQ_j<0$, so some $Q_j$, and hence at least two of them, are non-zero. Every horizon other than $S'$ satisfies $|Q_j|<1/(2\sqrt\Lambda)$, and therefore $|Q_{S'}|=\bigl|\sum_{\partial_jM\neq S'}Q_j\bigr|<(\ell-1)/(2\sqrt\Lambda)$.

If $M$ is not simply connected, then by Theorem~\ref{thm:D}(v) its universal covering $\tilde M$ is compact, with $\ell\,|\pi_1(M)|\ge4$ horizons. Since the horizons of $M$ are spheres, each horizon of $\tilde M$ is isometric to a horizon of $M$ and has the same charge, and it has an isometric copy among the other horizons of $\tilde M$; hence every horizon $\tilde S$ of $\tilde M$ satisfies $|\tilde S|\le|\partial\tilde M|-|\tilde S|$. As $\tilde M$ has at least four horizons, it is not a charged Nariai system, and the first part of the proof, applied to $\tilde M$, gives the strict bounds for every horizon of $M$.
\end{proof}

\begin{remark}\label{rem:largest}
The horizon of largest area is also controlled. If it is stable, the stability inequality gives $\Lambda|S|+16\pi^2Q_S^2/|S|\le4\pi$. If it is unstable, it has Morse index one by Theorem~\ref{thm:D}(iv), and the area--charge inequality for minimal surfaces of index one \cite[inequality~(5)]{CLS24}, which applies since $R_g=2\Lambda+2|E|^2$ and only involves the Jacobi operator of $S$ and the Gauss equation along $S$, gives $\Lambda|S|+16\pi^2Q_S^2/|S|\le12\pi$. Hence every horizon of a compact system satisfies $\Lambda|S|\le12\pi$ and, by the arithmetic-geometric mean inequality, $|Q_S|\le3/(2\sqrt\Lambda)$; for $\ell\ge5$ this improves the last bound of Theorem~\ref{thm:masscharge}.
\end{remark}

\begin{proof}[Proof of Theorem~\ref{thm:C}]
By Theorem~\ref{thm:D}, $\partial M\cong\Sph^2$, $H_1(M;\R)=0$, $\Lambda>0$, and $VE^\flat=d\psi$ with $\psi$ constant on the connected set $\partial M$; replacing $\psi$ by $\psi-\psi|_{\partial M}$ we may assume $\psi|_{\partial M}=0$. By Lemma~\ref{lem:basic}(iv) and the divergence theorem,
\[
\int_MV|E|^2\,dv_g=\int_M\dv(\psi E)\,dv_g=-\int_{\partial M}\psi\,\langle E,\nu\rangle\,d\sigma=0 ,
\]
so $E\equiv0$ because $V>0$ in $\inte M$. Then $R_g=2\Lambda$, and \eqref{eq:electrostatic} reduces to $\Hess V-(\Delta V)g-V\Ric=0$, i.e.\ $(M,g,V)$ is a static vacuum system in the sense of \cite{Amb17}.

We now prove (a)--(c). If $\partial M$ is unstable, Theorem~\ref{thm:D}(iv) shows that $M$ is a $3$-ball and that $\partial M$ has Morse index one. If $\partial M$ is stable and $M$ is simply connected, Theorem~\ref{thm:D}(v) shows that $M$ is a $3$-ball. In both cases $M$ is simply connected, and \cite[Theorem~D]{Amb17}, stated there for $R_g=6$, shows that the static vacuum system $(M,g,V)$ is either the de Sitter system or satisfies $\Lambda|\partial M|<4\pi$, a scale-invariant condition. Since the boundary of the hemisphere is unstable, this proves (a) and (b). In particular, if $M$ is not simply connected, then $\partial M$ is stable, and the proof of the inequality in Proposition~\ref{prop:rigidity}, which only uses stability, gives $\Lambda|\partial M|\le4\pi$.

Assume now that $M$ is not simply connected. By Theorem~\ref{thm:D}(v), $M$ is diffeomorphic to $\Sph^3/G$ with a ball removed, where $G=\pi_1(M)$ is finite and non-trivial. Let $M'\to M$ be the universal covering, of degree $d=|G|\ge2$. Since $\partial M$ is a sphere, its preimage consists of $d$ components isometric to $\partial M$, and $(M',g',V',0)$ is a compact electrostatic system with $d$ horizons of equal area. By Theorem~\ref{thm:D}(iii), each of them minimizes area in its homology class in $M'$. If $\Lambda|\partial M|=4\pi$, equality holds in Proposition~\ref{prop:rigidity} for a horizon of $M'$, so $M'$ is of product type; being compact with $d\ge2$ horizons and $E=0$, it is the Nariai system, so $d=2$. By the proof of Theorem~\ref{thm:B}, the only free orientation-preserving involution of the Nariai system preserving $V$ is $(s,p)\mapsto(-s,-p)$, so $M$ is the quotient of the Nariai system. Conversely, that quotient has $|\partial M|=4\pi/\Lambda$.

The assertion on projective planes is the case $\ell=1$ of Corollary~\ref{cor:rp2} below, whose proof does not use it. Since the quotient of the Nariai system contains a projective plane, equality in $\Lambda|\partial M|\le4\pi$ is impossible when $M$ contains none.
\end{proof}

The following corollary contains the inequality for projective planes in Theorem~\ref{thm:C}(c), for any number of horizons and all one-sided surfaces.

\begin{corollary}\label{cor:rp2}
Let $(M,g,V,E)$ be a compact electrostatic system with $V^{-1}(0)=\partial M$ consisting of $\ell\ge1$ horizons, and let $P\subset\inte M$ be a closed embedded one-sided surface. Then $|P|\ge\frac12|\partial M|$, with equality only if $(M,g,V,E)$ is the quotient of the Nariai system and $P$ is a projective plane, where it is attained. If $P$ is a projective plane, then $M$ is diffeomorphic to $\RP^3$ with $\ell$ open balls removed.
\end{corollary}

\begin{proof}
Since $P$ is one-sided, $M\setminus P$ is connected, and the mod $2$ intersection number with $P$ defines a connected double covering $M'\to M$, trivial over $M\setminus P$. Thus $M'\setminus\hat P$ consists of two copies of $M\setminus P$, where $\hat P$, the preimage of $P$, is the orientation double cover of $P$. Each copy contains one lift of each horizon, so $\hat P$ is homologous to the union $\partial_JM'$ of the lifts contained in one copy, and $|\partial_JM'|=|\partial M|=|\partial M'\setminus\partial_JM'|$. By Theorem~\ref{thm:D}(iii) applied to $M'$, $2|P|=|\hat P|\ge|\partial M|$, and equality forces $\hat P\subset\inte M'$ to be a minimizer, hence locally area-minimizing. Then Theorem~\ref{thm:B} applies to $P$, and $M$, being compact with non-empty boundary and not simply connected, is the quotient (viii) of the Nariai system; in it, $\hat P$ is a slice $\{s\}\times\Sph^2$, so $P$ is a projective plane. Conversely, in that quotient the image of $\{0\}\times\Sph^2$ is a projective plane of area $\frac12|\partial M|$. Finally, if $P$ is a projective plane, a regular neighborhood of $P$ is $\RP^3$ with a ball removed, so $M$ has a summand $\RP^3$; in particular $\pi_1(M)\neq1$, and by Theorem~\ref{thm:D}(v) and the uniqueness of the prime decomposition, $M$ is $\RP^3$ with $\ell$ balls removed.
\end{proof}

When $\pi_1(M)\neq1$, passing to the universal covering removes the condition on the areas in Theorem~\ref{thm:D}(iii) for spheres.

\begin{corollary}\label{cor:spheres}
Let $(M,g,V,E)$ be a compact electrostatic system with $V^{-1}(0)=\partial M$ consisting of $\ell\ge1$ horizons, and let $S\subset\inte M$ be an embedded sphere. Then $S$ separates $M$, and the closure of at least one component $W_S$ of $M\setminus S$ is a $3$-ball with balls removed.
\begin{itemize}
\item[(i)] If $M$ is simply connected, then $|S|\ge\min\{|\partial M\cap W|,|\partial M\cap W'|\}$, where $W,W'$ are the two components of $M\setminus S$.
\item[(ii)] If $M$ is not simply connected, $W_S$ is unique and $|S|\ge|\partial M\cap W_S|$.
\end{itemize}
If the right-hand side is positive, equality holds only if $(M,g,V,E)$ is a charged Nariai system in (i), and the quotient of the Nariai system in (ii). In particular, every embedded sphere in $\inte M$ which does not bound a ball has area at least $\min_j|\partial_jM|$, and, if $\pi_1(M)\neq1$, every embedded sphere separating all horizons from the summand $\Sph^3/G$ has area at least $|\partial M|$; the infimum of the areas of such spheres is $|\partial M|$.
\end{corollary}

\begin{proof}
By Theorem~\ref{thm:D}(v), $M$ is $\Sph^3/G$ with $\ell$ balls removed, so $H_1(M;\R)=0$ and $S$ separates $M$. Capping the horizons gives $\Sph^3/G$, which is irreducible by Alexander's theorem \cite{Ale24}, applied in the universal covering $\Sph^3$, so $S$ bounds a ball $B$ in $\Sph^3/G$, and $W_S=B\cap M$ is a ball with balls removed; if $G=1$, both components of $M\setminus S$ are of this type, again by Alexander's theorem. If $G\neq1$, the other component has fundamental group $G$, so $W_S$ is unique. In case (i), $S$ is homologous to $\partial M\cap W$, up to orientation, and Theorem~\ref{thm:D}(iii) gives the inequality; equality forces $S\subset\inte M$ to be a minimizer, so the system is a charged Nariai system. In case (ii), $\overline{W_S}$ is simply connected, so its preimage in the universal covering $\tilde M$ consists of $|G|\ge2$ disjoint copies, each bounded by a lift of $S$ and by lifts of the horizons in $W_S$, of total area $A=|\partial M\cap W_S|$. The other horizons of $\tilde M$ have total area $|G|\,|\partial M|-A\ge A$, so Theorem~\ref{thm:D}(iii), applied to $\tilde M$, gives $|S|\ge A$, with equality only if a lift of $S$ is a minimizer in $\inte\tilde M$, that is, only if $\tilde M$ is a charged Nariai system; then $\ell|G|=2$, and $M$ is the quotient of the Nariai system, as in the proof of Theorem~\ref{thm:C}. If $S$ does not bound a ball, every component of $M\setminus S$ which is a ball with balls removed contains a horizon; this gives the first bound in the last assertion, in both cases. If $S$ separates all horizons from the summand $\Sph^3/G$, then $\partial M\subset W_S$. Joining parallel copies of the horizons by thin tubes gives such spheres with area arbitrarily close to $|\partial M|$.
\end{proof}

\begin{proof}[Proof of Corollary~\ref{cor:gap}]
The three cases of Theorem~\ref{thm:C} exhaust all possibilities, and the de Sitter system has $\Lambda|\partial M|=12\pi$.
\end{proof}

\begin{proof}[Proof of Theorem~\ref{thm:E}]
Suppose first $\Lambda\le0$. Integrating $\Delta V=(|E|^2-\Lambda)V$ over the closed manifold $M$ (where $V>0$) gives $\int_M(|E|^2-\Lambda)V=0$ with a non-negative integrand; hence $|E|^2\equiv\Lambda$, so $\Lambda=0$ and $E\equiv0$. Then $V$ is harmonic, hence constant, and \eqref{eq:electrostatic} gives $V\Ric=\Hess V=0$, so $g$ is flat and, by Bieberbach's theorem \cite{Bie11}, $(M,g)$ is a finite quotient of a flat $3$-torus.

Suppose now $\Lambda>0$. If $H_2(M;\R)=0$, then $H^1(M;\R)\cong H_2(M;\R)=0$ by Poincar\'e duality, so $VE^\flat=d\psi$ by Lemma~\ref{lem:basic}(iv), and $\int_MV|E|^2=\int_M\dv(\psi E)=0$ forces $E\equiv0$; then $\Delta V=-\Lambda V$ gives $\Lambda\int_MV=0$, which is absurd. Therefore $H_2(M;\R)\neq0$, and there is $\alpha\in H_2(M;\Z)$ of infinite order. By Theorem~\ref{thm:homology}, $\alpha$ is represented by $\sum n_i[\Sigma_i]$ with $\Sigma_i$ closed, embedded, orientable and locally area-minimizing. By Theorem~\ref{thm:B} the system is one of (i)--(viii), and the only closed ones are those in (v).
\end{proof}

\section{Noncompact electrostatic systems}\label{sec:noncompact}

Unless stated otherwise, in this section $(M,g,V,E)$ is a noncompact electrostatic system with $V^{-1}(0)=\partial M$, mean-convex ends and exhaustion $\{\Omega_k\}$. Theorem~\ref{thm:F} assumes that the ends are strictly mean-convex; this is used only through Remark~\ref{rem:sharp}, to exclude systems of product type. Since $M$ is connected, we may assume that each $\Omega_k$ is connected, replacing it, for $k$ large, by its component containing $\partial M$ and a fixed point of $M$; these components still exhaust $M$. Then $\partial_+\Omega_k\neq\emptyset$, since otherwise $\Omega_k$ would be open and closed in $M$, hence equal to $M$, which is not compact. As $\partial M$ is totally geodesic, $\Omega_k$ is mean-convex.

\begin{lemma}\label{lem:noproduct}
Assume that $(M,g,V,E)$ is not of product type, and let $W\subset M$ be a compact mean-convex domain each of whose boundary components is either a component of $\partial M$ or contained in $\inte M$. Then every component of a mass minimizing integral cycle in $W$, as in Theorem~\ref{thm:homology}, is a component of $\partial W\cap\partial M$.
\end{lemma}

\begin{proof}
Let $\sum_in_i[\Sigma_i]$ be such a minimizer, and suppose that some $\Sigma_i$ is not a component of $\partial W\cap\partial M$, so that $\Sigma_i\subset\inte M$. By Theorem~\ref{thm:homology}, $\Sigma_i$ is closed, embedded and orientable, and either lies in $\inte W$ and is locally area-minimizing, or is a minimal component of $\partial W$. In both cases it has constant area along the flow: in the first by local minimality, and in the second for the flow $\Phi$ in the direction of the unit normal pointing into $W$, since for small $t>0$ the surface $\Phi(t,\Sigma_i)$ lies in $\inte W$, is homologous to $\Sigma_i$ in $W$ and is disjoint from the other $\Sigma_j$, so that replacing $n_i[\Sigma_i]$ by $n_i[\Phi(t,\Sigma_i)]$ gives a competitor and $|\Phi(t,\Sigma_i)|\ge|\Sigma_i|$. By Proposition~\ref{prop:globalsplit}, the system is of product type, contrary to the assumption.
\end{proof}

\begin{proof}[Proof of Theorem~\ref{thm:F}(i)]
By Remark~\ref{rem:sharp}, the system is not of product type, so Lemma~\ref{lem:noproduct} applies.

\emph{Stability.} Suppose that a component $\partial_1M$ is unstable, let $\lambda_1<0$ be the first eigenvalue of $L_{\partial_1M}$ and $u>0$ a first eigenfunction. For small $\varepsilon>0$, $\Sigma_\varepsilon=\{\exp_x(\varepsilon u(x)\nu(x)):x\in\partial_1M\}$ is an embedded surface in $\inte M$ bounding, together with $\partial_1M$, a collar $U_\varepsilon$. Its mean curvature $H_\varepsilon$ with respect to $\nu$ satisfies $\partial_\varepsilon H_\varepsilon|_{\varepsilon=0}=-L_{\partial_1M}u=\lambda_1u<0$, so $\vec H_\varepsilon$ is nowhere zero and points away from $U_\varepsilon$. Choose $k$ with $\overline{U_\varepsilon}\subset\Omega_k$ and let $W=\Omega_k\setminus(U_\varepsilon\cup\partial_1M)$. Then $W$ is compact, connected and mean-convex, with
\[
\partial W=\Sigma_\varepsilon\ \sqcup\ (\partial M\setminus\partial_1M)\ \sqcup\ \partial_+\Omega_k ,\qquad \partial_+\Omega_k\neq\emptyset .
\]
By Lemma~\ref{lem:boundaryclasses}, $\alpha=[\Sigma_\varepsilon]\in H_2(W;\Z)$ does not lie in the span of the classes of the components of $\partial M\setminus\partial_1M$. By Lemma~\ref{lem:noproduct}, however, a minimizer of $\alpha$ is supported on components of $\partial M\setminus\partial_1M$, a contradiction.

\emph{Strict outer-minimization.} Let $\Theta$ be a non-empty union of components of $\partial M$ and $\tilde\Sigma\neq\Theta$ a closed embedded orientable surface with $\tilde\Sigma-\Theta=\partial C$ for an integral $3$-chain $C$. Choose $k$ with $\tilde\Sigma\cup\operatorname{supp}C\subset\Omega_k$ and $\tilde\Sigma\cap\partial_+\Omega_k=\emptyset$, so that $[\tilde\Sigma]=[\Theta]$ in $H_2(\Omega_k;\Z)$. By Lemma~\ref{lem:noproduct}, a minimizer in this class has the form $\sum_jn_j[\partial_{i_j}M]$. As $\partial_+\Omega_k\neq\emptyset$, the classes of the components of $\partial M$ are linearly independent in $H_2(\Omega_k;\R)$ (Lemma~\ref{lem:boundaryclasses}), so the minimizer is $\Theta$ with multiplicity one and $|\Theta|\le|\tilde\Sigma|$. If equality holds, $\tilde\Sigma$ is itself a minimizer, so by Lemma~\ref{lem:noproduct} it is a union of components of $\partial M$, and hence $\tilde\Sigma=\Theta$, a contradiction.
\end{proof}

\begin{corollary}\label{cor:noncompactcharge}
Under the assumptions of Theorem~\ref{thm:F}(i), every component $S$ of $\partial M$ satisfies
\[
\Lambda|S|+\frac{16\pi^2Q_S^2}{|S|}<2\pi\chi(S),\qquad Q_S=\frac1{4\pi}\int_S\langle E,\nu\rangle .
\]
In particular: if $\Lambda\ge0$, every component of $\partial M$ is a sphere; if $\Lambda=0$, then $|S|>4\pi Q_S^2$; if $\Lambda>0$, then $4\Lambda Q_S^2<1$, $r_-(\Lambda,Q_S)<r_S<r_+(\Lambda,Q_S)$ and $m_-(\Lambda,Q_S)<m_H(S)<m_+(\Lambda,Q_S)$; and if $\Lambda<0$ and $S$ has genus $\mathrm g$, then
\[
|S|>\frac{2\pi}{|\Lambda|}\Bigl(\mathrm g-1+\sqrt{(\mathrm g-1)^2+4|\Lambda|Q_S^2}\Bigr).
\]
\end{corollary}

\begin{proof}
By Theorem~\ref{thm:F}(i), $|S|<|S'|$ for every closed embedded orientable surface $S'\neq S$ homologous to $S$, so Proposition~\ref{prop:rigidity} gives the inequality with $\le$, and equality would make the system of product type, which is excluded by Remark~\ref{rem:sharp}. If $\Lambda\ge0$, the left-hand side is non-negative, so $\chi(S)>0$ and $S$, being orientable, is a sphere. For $\Lambda=0$ the inequality reads $|S|>4\pi Q_S^2$. For $\Lambda>0$ it reads $\Lambda r_S^4-r_S^2+Q_S^2<0$, and the proof of Theorem~\ref{thm:masscharge} gives $4\Lambda Q_S^2<1$, $r_-<r_S<r_+$ and $m_-<m_H(S)<m_+$. For $\Lambda<0$, multiplying by $|S|$ gives $|\Lambda||S|^2-4\pi(\mathrm g-1)|S|-16\pi^2Q_S^2>0$, and the bound follows.
\end{proof}

For stable minimal or marginally trapped surfaces, inequalities of this type, with $\le$, are known; see \cite{Woo99} for the genus bounds when $\Lambda<0$ and \cite{DJR12,Sim12} for the charged case. The point of Corollary~\ref{cor:noncompactcharge} is that the horizons are automatically stable and that equality is excluded. For a connected horizon with $H_2(M,\partial M)=0$, the lower bound for $|S|$, as a non-strict inequality with equality only for a half-cylinder $[0,\infty)\times\partial M$, also follows from \cite[Theorem~3]{CM25}, proved by $\mu$-bubbles for complete initial data with $R_g\ge2\Lambda+2|E|^2$, $\Lambda+|E|^2$ uniformly positive and weakly mean-convex boundary, without assumptions on the ends. Corollary~\ref{cor:noncompactcharge} applies to every component of $\partial M$ and gives two-sided strict bounds; the half-cylinders in the equality case of \cite[Theorem~3]{CM25} that are electrostatic systems are the cold and ultracold black holes (ii) and (vi), whose ends are not strictly mean-convex.

For $\Lambda=0$, little is known about static electrovacuum solutions without asymptotic assumptions \cite{Rei24}. By Corollary~\ref{cor:noncompactcharge} and Remark~\ref{rem:sharp}, a noncompact system with $\Lambda=0$ and mean-convex ends which has a horizon of positive genus is of product type, and in a system of product type with $\Lambda=0$, a horizon and cross-sections of positive genus, $K=E_0^2\le0$ forces $E_0=0$ and $\mu=0$; hence it is the Rindler-type system (vi$'$) of Section~\ref{sec:models}, and $E\equiv0$. This is the boost black hole in Reiris' classification of metrically complete static vacuum black holes, according to which such a black hole is a boost, a Schwarzschild black hole, or of Myers--Korotkin--Nicolai type \cite{Rei18}.

\begin{proof}[Proof of Theorem~\ref{thm:F}(ii)]
As in the proof of (i), Lemma~\ref{lem:noproduct} applies. Let $\alpha\in H_2(M;\Z)$ and choose $k$ such that a cycle representing $\alpha$ lies in $\Omega_k$, with class $\alpha_k\in H_2(\Omega_k;\Z)$. By Lemma~\ref{lem:noproduct}, a minimizer of $\alpha_k$ has the form $\sum_jn_j[\partial_jM]$, so $\alpha=\sum_jn_j[\partial_jM]$ in $H_2(M;\Z)$.

Suppose that $M$ has at least two ends, and choose $k$ such that $M\setminus\Omega_k$ has two unbounded components $U_1$ and $U_2$. Since $M$ is connected, $Y_j=\partial\Omega_k\cap\overline{U_j}$ is a non-empty union of components of $\partial_+\Omega_k$. By Lemma~\ref{lem:boundaryclasses}, $[Y_1]$ is not a combination of the classes of the components of $\partial M$ in $H_2(\Omega_k;\Z)$, because in such a relation the components of $Y_1$ and of $Y_2$ would have the same coefficient. By Lemma~\ref{lem:noproduct}, however, a minimizer of $[Y_1]$ is supported on components of $\partial M$, a contradiction.

Let $\Sigma\subset\inte M$ be a closed embedded orientable surface. If $M\setminus\Sigma$ were connected, a closed curve $\gamma\subset\inte M$ meeting $\Sigma$ transversally in one point would have intersection number $\pm1$ with $[\Sigma]\in H_2(M;\Z)$, for the intersection pairing $H_2(M;\Z)\times H_1(M,\partial M;\Z)\to\Z$; but $[\Sigma]=\sum_jn_j[\partial_jM]$, and each $\partial_jM$ is homologous to a parallel copy disjoint from $\gamma$, a contradiction. Finally, since $\dv E=0$, the flux of $E$ through a closed surface depends only on its homology class.
\end{proof}

\begin{proof}[Proof of Theorem~\ref{thm:F}(iii)]
Let $\iota:\partial_1M\to M$ be the inclusion, suppose that $\ker(\iota_*)\neq1$, and let $c\subset\partial_1M$ be a loop which is essential in $\partial_1M$ and null-homotopic in $M$. A null-homotopy has compact image, so for $k$ large $c$ is null-homotopic in $\Omega_k$, and $\ker(\pi_1(\partial_1M)\to\pi_1(\Omega_k))\neq1$. By the loop theorem (Theorem~\ref{thm:loop}) there is a Jordan curve $\gamma\subset\partial_1M$, essential in $\partial_1M$, bounding an embedded disk in $\Omega_k$. By Theorem~\ref{thm:plateau}, $\gamma$ bounds a least area disk $D$ in $\Omega_k$ which is almost properly embedded or contained in $\partial\Omega_k$. In the latter case $D\subset\partial_1M$, and $\gamma$ would be null-homotopic in $\partial_1M$. Hence $D\cap\partial\Omega_k=\partial D=\gamma\subset\partial M$, so the compact set $D$ is disjoint from $\partial_+\Omega_k$. Disks bounded by $\gamma$ in a small neighborhood of $D$ therefore lie in $\Omega_k$, so $D$ has least area among them and satisfies the hypothesis of Theorem~\ref{thm:A}. By Theorem~\ref{thm:A}, $M$ is the de Sitter system, which is compact, a contradiction. The last assertion follows because $\pi_1$ of a closed orientable surface of positive genus would inject into $\pi_1(M)$.
\end{proof}

\begin{proof}[Proof of Theorem~\ref{thm:F}(iv)]
If $M$ is compact, set $\Omega_k=M$. Suppose first that an annulus $A\subset\Omega_k$, for some $k$, has least area among the annuli in $\Omega_k$ bounded by $\Gamma_1\cup\Gamma_2$ and is disjoint from $\partial_+\Omega_k$; if $A$ has least area among the annuli in $M$, this holds for $k$ large. By Theorem~\ref{thm:douglas}, $A$ is almost properly embedded or contained in $\partial\Omega_k$, and the latter is impossible, because a connected subset of $\partial\Omega_k$ lies in a single component while $\partial A$ meets two. Since $A$ is disjoint from $\partial_+\Omega_k$, the annuli bounded by $\Gamma_1\cup\Gamma_2$ in a small neighborhood of $A$ lie in $\Omega_k$, so $A$ satisfies the hypothesis of Theorem~\ref{thm:A}; but then, by the proof of Theorem~\ref{thm:A}, $A$ would be a hemisphere, which is impossible. This proves the first assertion.

If the strict Douglas condition held, we could choose an annulus $A_0$ with $|A_0|<\inf(|D_1|+|D_2|)$, the infimum being taken over disks $D_i\subset M$ bounded by $\Gamma_i$, and $k$ with $A_0\subset\Omega_k$ and $\Gamma_1,\Gamma_2$ null-homotopic in $\Omega_k$. Infima over $\Omega_k$ are not smaller than infima over $M$, so the strict Douglas condition holds in the mean-convex domain $\Omega_k$, and Theorem~\ref{thm:douglas} gives an annulus $A$ of least area among the annuli in $\Omega_k$ bounded by $\Gamma_1\cup\Gamma_2$. As above, $A$ is almost properly embedded, so $A\cap\partial\Omega_k=\partial A\subset\partial M$ and $A$ is disjoint from $\partial_+\Omega_k$, which is impossible by the first part.
\end{proof}

Theorem~\ref{thm:A} also gives lower bounds for the area of disks and annuli with boundary on the horizons.

\begin{corollary}\label{cor:disks}
Let $(M,g,V,E)$ be an electrostatic system with $V^{-1}(0)=\partial M$ and mean-convex ends, let $S$ be a component of $\partial M$, and let $\gamma\subset S$ be a smooth Jordan curve which bounds a disk in $S$; denote by $a_S(\gamma)$ the least area of such a disk. Then every immersed disk $D\subset M$ with $\partial D=\gamma$ satisfies $|D|\ge a_S(\gamma)$, with equality only if $D\subset S$ or $(M,g,V,E)$ is the de Sitter system, $\gamma$ is a great circle and $D$ is a totally geodesic hemisphere. If $\Gamma_1,\Gamma_2$ are such curves on distinct components $S_1,S_2$ of $\partial M$, every immersed annulus $A\subset M$ with $\partial A=\Gamma_1\cup\Gamma_2$ satisfies $|A|>a_{S_1}(\Gamma_1)+a_{S_2}(\Gamma_2)$, and the right-hand side is the infimum of the areas of such annuli.
\end{corollary}

\begin{proof}
Let $W=M$ if $M$ is compact, and otherwise $W=\Omega_k$ with $k$ so large that $D\subset\Omega_k$. Since $\gamma$ is null-homotopic in $W$, Theorem~\ref{thm:plateau} gives a least area disk $D_0$ in $W$ bounded by $\gamma$, which is almost properly embedded or contained in $\partial W$. In the first case, as in the proof of Theorem~\ref{thm:F}(iii), $D_0$ satisfies the hypothesis of Theorem~\ref{thm:A} and its proof, which show that the system is the de Sitter system and that $D_0$ is a totally geodesic hemisphere; then $\gamma$ is a great circle of $S$, which bounds two disks of area $|D_0|$ in $S$. In the second case $D_0\subset S$, so $|D_0|=a_S(\gamma)$. In both cases $|D|\ge|D_0|=a_S(\gamma)$, and if equality holds, $D$ is itself a least area disk in $W$, so it lies in $S$ or the first case occurs. For the annulus, the first part gives $\inf|D_i|=a_{S_i}(\Gamma_i)$ over disks $D_i$ bounded by $\Gamma_i$. Joining least area disks $D_i\subset S_i$ by a thin tube gives annuli bounded by $\Gamma_1\cup\Gamma_2$ with area arbitrarily close to $a_{S_1}(\Gamma_1)+a_{S_2}(\Gamma_2)$, and, since the strict Douglas condition fails (Theorem~\ref{thm:F}(iv), which holds for mean-convex ends by Remark~\ref{rem:sharp}), no annulus has smaller area. Thus $a_{S_1}(\Gamma_1)+a_{S_2}(\Gamma_2)$ is the infimum, and it is not attained, again by Theorem~\ref{thm:F}(iv).
\end{proof}

\section{Final remarks}\label{sec:related}

\emph{Several horizons.} Static solutions with $\Lambda>0$ and more than two horizons were speculated on by Boucher, Gibbons and Horowitz (see \cite[Section~1]{Amb17}), and Dias, Gibbons, Santos and Way \cite{DGSW23} found numerically, without a rigorous proof, static vacuum solutions with two black holes and a cosmological horizon, so $\ell=3$. For such systems, Theorems~\ref{thm:D} and~\ref{thm:masscharge}, Corollaries~\ref{cor:spheres} and~\ref{cor:disks} and Theorem~\ref{thm:F}(iv) give necessary conditions, and $\inte M$ contains no closed locally area-minimizing surface, since the compact models of Theorem~\ref{thm:B} have at most two horizons.

\emph{Two horizons.} The uniqueness of the Schwarzschild--de Sitter and Reissner--Nordstr\"om--de Sitter systems is open in general; for $E=0$ it is known under assumptions on the level sets of $V$ \cite{BCM19}.

\begin{remark}[Near-horizon geometries]\label{rem:nearhorizon}
The spacetimes $(\R\times M,-V^2dt^2+g)$ associated with the models (i)--(vii) of Section~\ref{sec:models} are products of a two-dimensional Lorentzian space form with a round sphere, since $-V(s)^2dt^2+ds^2$ has constant curvature when $V''=\mu V$, $\mu=E_0^2-\Lambda$: one obtains $AdS_2\times\Sph^2$ for (i), (ii) and (iii), in global, Rindler and Poincar\'e coordinates, $\R^{1,1}\times\Sph^2$, the Pleba\'nski--Hacyan spacetime \cite{PH79}, for (iv) and (vi), and $dS_2\times\Sph^2$ for (vii). Theorem~\ref{thm:B} characterizes the static slices of these spacetimes variationally.
\end{remark}

\emph{Open problems.}
\begin{enumerate}[label=(\arabic*),leftmargin=*,itemsep=4pt,topsep=4pt]
\item The Boucher--Gibbons--Horowitz conjecture asserts that the de Sitter system is the only compact simply connected static vacuum system with connected horizon. By Theorem~\ref{thm:C}, every compact electrostatic system with connected horizon is a static vacuum system, and it is simply connected exactly in cases (a) and (b). Thus the conjecture says that the de Sitter system is the only system in cases (a) and (b). By Corollary~\ref{cor:gap}, a counterexample would satisfy $\Lambda|\partial M|<4\pi$, a third of the value for the de Sitter system. Does a compact static vacuum system on the $3$-ball with connected horizon and $\Lambda|\partial M|<4\pi$ exist?
\item Characterize the systems in which $|\partial_JM|=|\partial M\setminus\partial_JM|$ for some union $\partial_JM$ of horizons, the equality case of Theorem~\ref{thm:D}(iii); then both $\partial_JM$ and $\partial M\setminus\partial_JM$ minimize area in their homology class. For $\ell=2$: are the charged Nariai systems the only compact systems with two horizons of equal area? For $\ell\ge3$, no compact system is known to us.
\item By Theorem~\ref{thm:D}(iii), two horizons of equal area are both stable. Is the converse true, that is, does every compact system with two stable horizons satisfy $|\partial_1M|=|\partial_2M|$? This holds for the Reissner--Nordstr\"om--de Sitter systems, whose cosmological horizon is unstable \cite[Section~3.8]{CLS24}, and for the charged Nariai systems, whose horizons have equal area.
\item By Theorem~\ref{thm:D}(v), a compact system with horizons is $\Sph^3/G$ with balls removed. Which groups $G$ occur? The known examples have $G=1$ or, for the quotient of the Nariai system, $G=\Z_2$, which occurs exactly when $M$ contains a projective plane (Corollary~\ref{cor:rp2}).
\end{enumerate}

\appendix
\section{Local rigidity at an area-minimizing horizon}\label{app:rigidity}

The following result is due to Cruz and Mendes \cite[Propositions~4 and~6]{CM25}, except for the last assertion, which follows from Theorem~\ref{thm:B}. It is the boundary version of the local rigidity theorems of Bray, Brendle and Neves \cite{BBN10}, Cai and Galloway \cite{CG00}, Nunes \cite{Nun13} and M\'aximo and Nunes \cite{MN13}, in the charged setting of \cite{GM25,Men25} (see also \cite{BBB23}).

\begin{proposition}\label{prop:rigidity}
Let $(M,g,V,E)$ be an electrostatic system with $V^{-1}(0)=\partial M$, and let $S$ be a compact component of $\partial M$ such that $|S|\le|S'|$ for every closed embedded surface $S'\subset M$ homologous to $S$ in a neighborhood of $S$. Let $Q=\frac1{4\pi}\int_S\langle E,\nu\rangle$. Then $S$ is stable and
\[
\Lambda|S|+\frac{16\pi^2Q^2}{|S|}\le2\pi\chi(S).
\]
If equality holds, then $(M,g,V,E)$ is of product type.
\end{proposition}

\begin{proof}
Since $S$ minimizes area among nearby homologous surfaces, it is stable. Taking $\varphi=1$ in the stability inequality and using that $S$ is totally geodesic, we get $0\le-\int_S\Ric(\nu,\nu)$. By the Gauss equation and Lemma~\ref{lem:basic}(i), $\Ric(\nu,\nu)=\Lambda+|E|^2-K$ on $S$, where $K$ is its Gaussian curvature, and, since $E$ is normal to $S$, the Cauchy--Schwarz inequality gives $\int_S|E|^2\ge16\pi^2Q^2/|S|$; hence
\[
\Lambda|S|+\frac{16\pi^2Q^2}{|S|}\le\Lambda|S|+\int_S|E|^2\le\int_SK=2\pi\chi(S),
\]
as in \cite[Proposition~4]{CM25}. Assume that equality holds. Since $R_g=2\Lambda+2|E|^2$, $\dv E=0$ and $S$ has least area among nearby homologous surfaces, \cite[Proposition~6]{CM25} gives a neighborhood of $S$ isometric to $([0,\delta)\times S,dt^2+g_0)$, where $g_0$ has constant Gaussian curvature and $E$ is normal to the slices. For $t_0\in(0,\delta)$, the slice $\{t_0\}\times S\subset\inte M$ is locally area-minimizing, since the projection $P$ onto $S$ does not increase area and $P\circ F(1,\cdot)$ is a diffeomorphism onto $S$ for every deformation $F$ of the slice of small size. By Theorem~\ref{thm:B}, the system is of product type.
\end{proof}

\section*{Statements and Declarations}
\emph{Funding.} The author did not receive support from any organization for the submitted work.

\emph{Competing interests.} The author has no competing interests to declare.

\emph{Data availability.} Data sharing is not applicable, as no datasets were generated or analyzed.

\emph{Use of AI tools.} An AI language model (Claude, Anthropic) was used to assist in formalizing the proofs, in finding applications of Theorems~\ref{thm:A} and~\ref{thm:B}, in suggesting arguments and references, in drafting and editing the text, and in checking the literature. The author verified all mathematical content and references and takes full responsibility for them.

\providecommand{\bysame}{\leavevmode\hbox to3em{\hrulefill}\thinspace}
\providecommand{\MR}{\relax\ifhmode\unskip\space\fi MR }
\providecommand{\MRhref}[2]{%
  \href{http://www.ams.org/mathscinet-getitem?mr=#1}{#2}
}
\providecommand{\href}[2]{#2}

\end{document}